\input amstex 
\documentstyle{amsppt} 
\hsize 13cm
\vsize 16.05cm
\magnification=\magstep1
\NoBlackBoxes
\def\nmb#1#2{#2}         
\def\cit#1#2{\ifx#1!\cite{#2}\else#2\fi} 
\def\totoc{}             
\def\idx{}               
\def\ign#1{}             

\define\>{\rightarrow}
\define\<{\leftarrow}
\define\[{\lbrack}
\define\]{\rbrack}
\redefine\o{\circ}

\define\al{\alpha}
\define\be{\beta}
\define\ga{\gamma}
\define\de{\delta}
\define\ep{\varepsilon}

\define\th{\theta}
\define\io{\iota}

\define\la{\lambda}

\define\ta{\tau}
\define\ph{\varphi}

\define\ps{\psi}
\define\om{\omega}

\define\De{\Delta}

\define\Si{\Sigma}

\redefine\i{^{-1}}
\define\row#1#2#3{#1_#2,\ldots,#1_#3}
\define\x{\times}
\define\ev{\operatorname{ev}}
\def\today{\ifcase\month\or
 January\or February\or March\or April\or May\or June\or
 July\or August\or September\or October\or November\or December\fi
 \space\number\day, \number\year}
\topmatter
\title The Convenient Setting for Real Analytic Mappings \endtitle
\author Andreas Kriegl \qquad {\eightpoint\rm and} \qquad
        Peter W. Michor  \endauthor
\affil  {\eightpoint Universit\"at Wien,\\
Wien Austria} \endaffil
\date{August 21, 1989}\enddate
\address
Institut f\"ur Mathematik, Universit\"at Wien,
Strudlhofgasse 4, A-1090 Wien, Austria
\endaddress
\abstract{We present here "the" cartesian closed theory for real
analytic mappings. It is based on the concept of real analytic curves in
locally convex vector spaces. A mapping is 
real analytic, if it maps smooth curves to smooth curves and
real analytic curves to real analytic curves. Under mild
completeness conditions 
the second requirement can be replaced by: real analytic along affine
lines. Enclosed and necessary is a careful study of locally
convex topologies on spaces of real analytic mappings.

As an application we also present the theory of manifolds of real
analytic mappings: the group of real analytic diffeomorphisms of
a compact real analytic manifold is a real analytic Lie group.}
\endabstract
\endtopmatter
\document
\block
{\eightpoint
\heading Contents \endheading
\noindent 0. Introduction 
	\leaders \hbox to 1em{\hss .\hss }\hfill {\eightrm 2}\par 
\noindent 1. Real analytic curves 
	\leaders \hbox to 1em{\hss .\hss }\hfill {\eightrm 3}\par 
\noindent 2. Real analytic mappings 
	\leaders \hbox to 1em{\hss .\hss }\hfill {\eightrm 11}\par 
\noindent 3. Function spaces in finite dimensions 
	\leaders \hbox to 1em{\hss .\hss }\hfill {\eightrm 18}\par 
\noindent 4. A uniform boundedness principle 
	\leaders \hbox to 1em{\hss .\hss }\hfill {\eightrm 26}\par 
\noindent 5. Cartesian closedness 
	\leaders \hbox to 1em{\hss .\hss }\hfill {\eightrm 30}\par 
\noindent 6. Consequences of cartesian closedness 
	\leaders \hbox to 1em{\hss .\hss }\hfill {\eightrm 36}\par 
\noindent 7. Spaces of sections of vector bundles 
	\leaders \hbox to 1em{\hss .\hss }\hfill {\eightrm 41}\par 
\noindent 8. Manifolds of analytic mappings 
	\leaders \hbox to 1em{\hss .\hss }\hfill {\eightrm 47}\par 
}\endblock
\heading\totoc 0. Introduction \endheading

We always wanted to know whether the group of real analytic
diffeomorphisms of a real analytic manifold is itself a real
analytic manifold in some sense. The paper \cit!{16}
contains the theorem, that this group for a compact real
analytic manifold is a smooth Lie group modeled on locally
convex vector spaces. (The proof, however, contains a gap, which
goes back to Smale in \cit!{1}: in canonical charts,
no partial mapping of the composition is linear off 0).
The construction there relies on ad hoc 
descriptions of the topology on the space of real analytic
functions. Also the literature dealing with the duals of these
spaces like \cit!{9} does not really try to
describe the topologies on spaces of real analytic functions.
There are, however, some older papers on this subject, see 
\cit!{29}, \cit!{25}, 
\cit!{26}, \cit!{27}, 
\cit!{31}, \cit!{12},  
\cit!{8}.

For some other instances where real analytic mappings in
infinite dimensions make their appearance, see the survey
article \cit!{28}.

In this article, we present a careful study of real analytic
mappings in infinite (and finite) dimensions combined with a
thorough treatment of locally convex topologies on spaces of
real analytic functions. From the beginning our aim is cartesian
closedness: a mapping $f: E\x F\to G$ should be real analytic if
and only if the canonically associated mapping $\check f:E\to
C^\om(F,G)$ is it. Very simple examples, see \nmb!{1.1}, show
that real analytic in the sense of having a locally
converging Taylor series is too restrictive. 

The right notion
turns out to be scalarwise real analytic: A curve in a locally
convex space is 
called (scalarwise) real analytic if and only if composed with
each continuous linear 
functional it gives a real analytic function. Later we show,
that the space of real analytic curves does not depend on the topology,
only on the bornology described by the dual.

A mapping will be called real analytic if it maps smooth
curves to smooth curves and real analytic curves to real analytic curves. 
This definition is in spirit very near to the original ideas of
variational calculus and it leads to a simple and powerful theory.
We will show the surprising result, that under some mild
completeness conditions (i.e. for convenient vector spaces), the
second condition can be replaced by: the mapping should be real analytic
along affine lines, see \nmb!{2.7}. This is a version of
Hartogs' theorem, which 
for Banach spaces is due to \cit!{2}. 

It is a very satisfying result, that the right realm of spaces
for real analytic analysis is the category of convenient vector
spaces, which is also the good setting in infinite dimensions for
smooth analysis, see \cit!{5}, and for
holomorphic analysis, see \cit!{15}.

The power of the cartesian closed calculus for real analytic
mappings developed here is seen in \cit!{21}, where it
is used to construct, for any unitary representation of any Lie
group, a real analytic moment mapping from the space of
analytic vectors into the dual of the Lie algebra.

We do not give any hard implicit function theorem in this paper,
because our setting is too weak to obtain one ---
but we do not think that this is a disadvantage. Let us make a
programmatic statement here:

An eminent mathematician once said, that for infinite
dimensional calculus each serious application needs its own foundation. 
By a serious application one obviously means some application of
a hard inverse function theorem. These theorems can be proved,
if by assuming enough a priori estimates one creates enough
Banach space situation for some modified iteration procedure to
converge. Many authors try to build their platonic idea of an a
priori estimate into 
their differential calculus. We think that this makes the calculus
inapplicable and hides the origin of the a priori estimates. We
believe, that the calculus itself should be as easy to use as
possible, and that all further assumptions (which most often
come from ellipticity of some nonlinear partial differential
equation of geometric origin) should be treated separately, in a
setting depending on the specific problem. We are sure that in
this sense 
the setting presented here (and the setting in \cit!{5}) is universally
usable for most applications.

The later parts of this paper are devoted to the study of
manifolds of real analytic mappings. We show indeed, that the
set of real analytic mappings from a compact manifold to another
one is a real analytic manifold, that composition is real
analytic and that the group of real analytic diffeomorphisms is
a real analytic Lie group. The exponential mapping of it
(integration of vector fields) is real analytic, but as in the
smooth case it is still not
surjective on any neighborhood of the identity.
We would like to stress the fact that the group of smooth
diffeomorphisms of a manifold is a smooth but {\it not} a real
analytic Lie group. We also show that the space of smooth
mappings between real analytic manifolds is a real analytic
manifold, but the composition is only smooth.

Throughout this paper our basic guiding line is the cartesian
closed calculus for smooth mappings as exposed in
\cit!{5}. The reader is assumed to be
familiar with at least the rudiments of it; but section \nmb!{1}
contains a short summary of the essential parts.

We want to thank Janusz Grabowski for  hints and discussions. This
should have been a joint work with him, but distance prevented it.

\heading\totoc \nmb0{1}. Real analytic curves
\endheading

\subheading{\nmb.{1.1}} As for smoothness and holomorphy we
would like to obtain cartesian closedness for real analytic
mappings. Thus one should have at least the following:

$f:\Bbb R^2 \to  \Bbb R$ is real analytic in the classical sense
if and only if $f\spcheck :\Bbb R \to  C^\om(\Bbb R,\Bbb R)$ is real
analytic in some appropriate sense.

The following example shows that there are some subtleties involved.

\proclaim{Example} The mapping 
$$f:\Bbb R^2\ni(s,t)\mapsto \frac1{(st)^2+1}\in\Bbb R$$
is real analytic, whereas there is no reasonable
topology on $C^\om(\Bbb R,\Bbb R)$, such that the mapping  
$f\spcheck :\Bbb R \to  C^\om(\Bbb R,\Bbb R)$ is locally given by
its convergent Taylor series.
\endproclaim
\demo{Proof} For a topology on $C^\om(\Bbb R,\Bbb R)$ to be
reasonable we require only that all evaluations 
$\operatorname{ev}_t: C^\om(\Bbb R,\Bbb R) \to  \Bbb R$ are
bounded linear functionals. 
Now suppose that $f\spcheck(s) = \sum_{k=0}^\infty f_ks^k$ converges
in $C^\om(\Bbb R,\Bbb R)$ for small $t$, where $f_k \in
C^\om(\Bbb R,\Bbb R)$. Then the series converges even
bornologically, see \nmb!{1.7}, so $f(s,t) =
\operatorname{ev}_t(f\spcheck(s)) = 
\sum f_k(t)\,s^k$ 
for all $t$ and small $s$. On the other hand 
$f(s,t) = \sum_{k=0}^\infty (-1)^k(st)^{2k}$ for $|s|<1/{|t|}$.
So for all $t$ we have $f_k(t)=(-1)^mt^k$ for $k=2m$, and $0$
otherwise, since for fixed $t$ we have a real analytic function
in one variable. Moreover, the series
$\left(\sum f_kz^k\right)(t) = \sum (-1)^kt^{2k}z^{2k}$
has to converge in $C^\om(\Bbb R,\Bbb R)\otimes\Bbb C$ for 
$|z|\leq \de$ and all $t$, see \nmb!{1.7}. This is not the case:
use $z=\sqrt{-1}\,\de$, 
$t=1/\de$. 
\qed\enddemo

\subheading{\nmb.{1.2}} There is, however, 
another notion of real analytic curves. \newline 
{\bf Example.} Let $f:\Bbb R\to \Bbb R$
be a real analytic function with finite radius of convergence at
$0$. Now consider the curve $c:\Bbb R \to  \Bbb R^{\Bbb N}$
defined by $c(t):= (f(k\cdot t))_{k\in \Bbb N}$. Clearly the
composite of $c$ with any continuous linear functional is real
analytic, since these functionals depend only on finitely many
coordinates. But the Taylor series of $c$ at $0$ has radius of
convergence $0$, since the radii of the coordinate functions go
to $0$. For an even more natural example see \nmb!{5.2}.

The natural setting for this notion of real analyticity is that of
dual pairs:

\subheading{Definition (Real analytic curves)} Let a \idx{\it dual pair}
$(E,E')$ be a real 
vector space $E$ with prescribed point separating dual $E'$. 
A curve $c: \Bbb R\to  E$ is called 
\idx{\it real analytic} if $\la \o c:\Bbb R \to  \Bbb R$ is real
analytic for all $\la  
\in E'$.

A subset $B\subseteq E$ is called \idx{\it bounded} if $\la(B)$ is bounded in 
$\Bbb R$ for all $\la \in E'$. The set of bounded subset of $E$ will be 
called the \idx{\it bornology} of $E$ (generated by $E'$).

The dual pair $(E,E')$ is called \idx{\it complete} if the bornology on $E$ 
is complete, 
i.e. for every bounded set $B$ there exists a bounded absolutely convex 
set $A\supseteq B$ such that the normed space $E_A$ generated by
$A$, see \cit!{10}, 8.3 or 
\cit!{5}, 2.1.15, is
complete.

Let $\ta$ be a topology on $E$, which is compatible with the bornology 
generated by $E'$, 
i.e. has as von Neumann bornology exactly this bornology.
Then a curve $c:\Bbb R \to  (E,\ta)$ will be called 
\idx{\it topologically real analytic} if it is locally given by a 
power series converging with respect to $\ta$.

A curve $c: \Bbb R\to  E$ will be called 
\idx{\it bornologically real analytic} 
if it factors locally over a topologically real analytic curve into 
$E_B$ for some bounded absolutely convex set $B\subseteq E$.

\subheading{\nmb.{1.3} Review of the smooth and holomorphic setting}
We will make use of the cartesian closedness of smooth maps between 
convenient vector spaces \cit!{14} and that of holomorphic 
maps between such spaces \cit!{15}.
Let us recall some facts from those theories.

First the smooth theory, where we refer to \cit!{5}.
\idx{\it Separated preconvenient vector spaces} can be defined as those dual 
pairs $(E,E')$ for which $E'$ consists exactly of the linear
functionals which are boun\-ded with respect to the bornology on
$E$ generated by $E'$. 
 To each dual pair $(E,E')$ one can
naturally associate a preconvenient vector space $(E,E^b)$,
where $E^b$ denotes the space of linear functionals which are
boun\-ded for the bornology generated by $E'$.  The space
$(E,E^b)$ is the dual pair with the finest structure, which has
as underlying space $E$ and which has the same bornology.
On every dual pair there is a \idx{\it natural locally convex topology}, namely
the Mackey topology associated with $E'$. The associated
bornological topology given by the absolutely convex bornivorous
subsets of $E$ is the natural topology of $(E,E^b)$.
A curve $c:\Bbb R \to  E$ is called \idx{\it smooth} if $\la\o c:\Bbb R \to 
\Bbb R$ is smooth. If $(E,E')$ is complete and $\ta$ is any topology on 
$E$ that is compatible with the bornology, then $c$ is smooth if and only if $c$ has 
derivatives of arbitrary order with respect to $\ta$ or, equivalently,
for every $k$ the curve $c$ factors locally as a $\Cal Lip^k$-mapping over 
$E_B$ for some bounded absolutely convex set $B\subseteq E$.

A \idx{\it convenient vector space} or \idx{\it convenient dual pair} is a
separated preconvenient vector space 
$(E,E')$, which is complete, so that $E'=E^b$ and the natural
topology is bornological. Since the completeness condition
depends only on the bornology, $(E,E')$ is complete if and only
if $(E,E^b)$ is convenient.

A set $U\subseteq E$ is called \idx{\it $c^\infty$-open} if the inverse
image $c^{-1}(U)\subseteq\Bbb R$ is open for every smooth
curve $c$ or, equivalently, the intersection $U_B:=U\cap E_B$ is open
in the normed space $E_B$ for every bounded absolutely convex set 
$B\subseteq E$. If $E$ is a metrizable or a Silva locally convex 
space and $E'$ its topological dual then its topology  coincides with 
the $c^\infty$-topology.

A mapping $f:U\to F$ into another dual pair
$(F,F')$ is called \idx{\it smooth (or $C^\infty$)} if $f\o c$
is a smooth curve for every smooth curve $c$ having values in $U$.
For Banach or even Fr\'echet spaces this notion coincides with
the classically considered notions.  The \idx{\it space of smooth
mappings} from U to F will be denoted by $C^\infty(U,F)$. On
$C^\infty(\Bbb R,\Bbb R)$ we consider the Fr\'echet topology of
uniform convergence on compact subsets of all derivatives
separately. On $C^\infty(U,F)$ one considers the dual induced by
the family of mappings
$C^\infty(c,\la):C^\infty(U,F)\to C^\infty(\Bbb R,\Bbb R)$ for
$c\in C^\infty(\Bbb R,U)$ and $\la\in F'$. 
This makes $C^\infty(U,F)$ into a complete
dual pair provided F is complete, and so one can pass to the
associated convenient vector space. If $E$ and $F$ are finite 
dimensional the bornological topology of $C^\infty(U,F)$ is the
usual topology of 
uniform convergence on compact subsets of $U$ of all 
derivatives separately. For this space the following
\idx{\it exponential law} is valid: For every $c^\infty$-open set $V$ of a
convenient vector space a mapping $f:V\times U\to F$ is
smooth if and only if the associated mapping $\check
f:V\to C^\infty(U,F)$ is a well defined smooth map.

Now the holomorphic theory developed in \cit!{15}.
Let $\Bbb D$ denote the open unit disk $\{z\in\Bbb C:|z| < 1\}$
in $\Bbb C$.  For a complex dual pair $(E,E')$ a map
$c:\Bbb D\to E$ is called a \idx{\it holomorphic curve} if $\la\o c:\Bbb
D\to \Bbb C$ is a holomorphic function for every $\la \in E'$. If $(E,E')$ is 
complete and $\ta$ is any topology on $E$ that is compatible with the 
bornology, then $c$ is holomorphic if and only if $c$ is complex differentiable 
with respect to $\ta$ or, equivalently,
the mapping $c$ factors locally as a holomorphic curve over 
$E_B$ for some bounded absolutely convex set $B\subseteq E$.
A mapping $f:U\to F$ between
complete complex dual pairs is called \idx{\it holomorphic} if $f\o c:\Bbb
D\to  F$ is a holomorphic curve for every holomorphic curve $c$
having values in $U$. This is true if and only if it is a smooth
mapping for the 
associated real vector spaces and the derivative at every point
in $U$ is $\Bbb C$-linear. For Banach or even Fr\'echet spaces
this notion coincides with classically considered notions.
Let $\Cal H(U,F)$ denote the \idx{\it vector space of holomorphic maps}
from $U$ to $F$. Then $\Cal H(U,F)$ is a closed subspace of
$C^\infty(U,F)$, since $f\mapsto f'(x)(v)$ is continuous on the
latter space. So one equips $\Cal H(U,F)$ with the convenient
vector space structure induced from $C^\infty(U,F)$.
If $E$ is finite dimensional, then the bornological topology on
$\Cal H(U,F)$ is the topology of uniform convergence on compact
subsets of $U$, see \nmb!{3.2}.
For this space one has again an \idx{\it exponential law}: For every
$c^\infty$-open subset $V$ of a complex convenient vector space
a mapping $f:V\times U\to F$ is holomorphic if and only if the
associated mapping $\check f:V\to \Cal H(U,F)$ is a well defined
holomorphic map.  This is a slight generalization of
\cit!{15}, 2.14, with the same proof as given there.

\proclaim{\nmb.{1.4}. Lemma} For a formal
power series $\sum_{k\geq0}a_kt^k$ with real coefficients the
following conditions are equivalent. 
\roster
\item The series has positive radius of convergence.
\item $\sum a_kr_k$ converges absolutely for all sequences $(r_k)$ with
        $r_k\,t^k\to 0$ for all $t>0$.
\item The sequence $(a_kr_k)$ is bounded for all $(r_k)$ with
        $r_k\,t^k\to 0$ for all $t>0$.
\item For each sequence $(r_k)$ satisfying $r_k>0$, 
        $r_kr_\ell\geq r_{k+\ell}$, 
        and $r_k\,t^k\to 0$ for all $t>0$ there exists an
        $\ep>0$ such that $(a_k\,r_k\,\ep^k)$ is bounded.
\endroster
\endproclaim

This bornological description of real analytic curves will be
rather important for the theory presented here, since condition 
\therosteritem3 and \therosteritem4 are linear
conditions on the coefficients of a formal power series
enforcing local convergence.

\demo{Proof} \therosteritem1 $\Rightarrow$ \therosteritem2. 
$\sum a_kr_k = \sum (a_kt^k)(r_kt^{-k})$ converges absolutely
for some small $t$.

\therosteritem2 $\Rightarrow$ \therosteritem3
$\Rightarrow$ \therosteritem4 is clear.

\therosteritem4 $\Rightarrow$ \therosteritem1. If the series
has radius of convergence 0, then we have 
$\sum_k|a_k|\,(\frac1{n^2})^k = \infty$ for all $n$. 
There are $k_n\nearrow\infty$ with
$$\sum_{k=k_{n-1}}^{k_n-1} |a_k|\,(\tfrac1{n^2})^k \geq 1.$$
We put $r_k := (\frac1n)^k$ for $k_{n-1}\leq k < k_n$, then
$\sum_k|a_k|\,r_k(\frac1n)^k = \infty$ for all $n$, so
$(a_kr_k(\tfrac1{2n})^k)_k$ is not bounded for any $n$, 
but $r_k\,t^k$, which equals $(\frac tn)^k$ for
$k_{n-1}\leq k < k_n$, converges to 0 for all $t>0$, and the
sequence $(r_k)$ is subadditive as required.
\qed\enddemo

\proclaim{\nmb.{1.5}. Theorem (Description of real
analytic curves)}
Let $(E,E')$ be a complete dual pair. A
curve $c:\Bbb R\to E$ is real analytic if and only if $c$ is 
smooth, and for each sequence $(r_k)$ with $r_k\,t^k\to
0$ for all $t>0$, and each compact set $K$ in $\Bbb R$, 
the set 
$$\left\{\frac1{k!}\,c^{(k)}(a)\,r_k: a\in K, k\in \Bbb N\right\}$$
is bounded, or equivalently the set corresponding to \nmb!{1.4}.(4) is 
bounded.
\endproclaim
\demo{Proof} Since both conditions can be tested by applying 
$\la\in E'$ and we have $(\la\o c)^{(k)}(a) = \la(c^{(k)}(a))$ we may 
assume that $E=\Bbb R$.

($\Rightarrow$). Clearly $c$ is smooth.\newline
{\bf Claim.} There exist $M,\rho>0$ with $|\frac1{k!}c^{(k)}(a)|<M\rho^k$ 
for all $k\in\Bbb N$ and $a\in K$.

This will give that
$|\frac1{k!}\,c^{(k)}(a)\,r_k\ep^k|\leq Mr_k(\ep\rho)^k$  which is 
bounded since $r_k(\ep\rho)^k\to 0$, as required.

To show the claim we argue as follows.
Since the Taylor series of $c$ converges at $a$ there are constants 
$M_a,\rho_a$ satisfying the claimed inequality for fixed $a$. 
An elementary computation shows that for all $a'$ with 
$|a-a'|\leq\tfrac1{2\rho_a}$ we have
$$\left|\frac{c^{(k)}(a')}{k!}\right| \leq M_a\rho_a^k\;\frac1{k!}
\left.\frac{\partial^k}{\partial t^k}\right|_{t=\tfrac12}\frac1{1-t}\;,$$
hence the condition is satisfied for all those $a'$ 
with some new constants $M_a',\rho_a'$.
Since $K$ is compact the claim follows. 

($\Leftarrow$). 
Let 
$$a_k:=\sup_{a\in K}\left|\frac1{k!}\,c^{(k)}(a)\right|.$$
Using \nmb!{1.4} (4$\Leftarrow$1) 
these are the coefficients of a power series with positive 
radius $\rho$ of convergence. Hence the remainder of the Taylor series 
goes locally to zero.
\qed\enddemo

Although topological real analyticity is a strictly stronger than real 
analyticity, cf. \nmb!{1.2}, sometimes the converse is true as the 
following slight generalization of \cit!{2}, Lemma 7.1
shows.

\proclaim{\nmb.{1.6}. Theorem}  Let $(E,E')$ be a complete dual pair and 
assume that a Baire vector space topology on $E'$ exists for
which the point evaluations $\ev_x$ for $x\in E$ are continuous.
Then any real analytic curve $c:\Bbb R \to  E$ is locally  
given by its Mackey convergent Taylor series, and hence is 
bornologically real analytic and topologically real analytic for every 
locally convex topology compatible with the bornology.
\endproclaim
\demo{Proof} Since $c$ is real analytic, it is smooth and all
derivatives exist in $E$, since $(E,E')$ is complete, by
\nmb!{1.3}.

Let us fix $t_0 \in \Bbb R$, let $a_n := 
\frac 1{n!}c^{(n)}(t_0)$. 
It suffices to find some $r>0$ for which 
$\{r^n a_n: n\in \Bbb N_0\}$ is bounded; because then $\sum
t^n a_n$ is Mackey-convergent for $|t|<r$, and its limit is
$c(t_0+t)$ since we can test this with functionals.

Consider the sets $A_r := \{\la\in E': |\la(a_n)| \leq r^n
\text{ for all }n\in \Bbb N\}$. These $A_r$ are closed in the
Baire topology, since the point evaluations at $a_n$ are continuous.
Since $c$ is real analytic, $\bigcup_{r>0} A_r = E'$, and by the Baire
property there is an $r>0$ such that the interior $U$ of $A_r$
is not empty. Let $\la_0\in U$, then for all $\la$ in the open
neighborhood $U-\la_0$ of $0$ we have 
$|\la(a_n)|\leq |(\la+\la_0)(a_n)|+|\la_0(a_n)|\leq 2r^n$.
The set $U-\la_0$ is absorbing, thus for every $\la \in E'$
some multiple $\ep\la$ is in $U-\la_0$ and so 
$\la(a_n)\leq \frac2\ep r^n$ as required.
\qed\enddemo

\proclaim{\nmb.{1.7}. Lemma} Let $(E,E')$ be a complete dual pair, $\ta$ a 
topology on $E$ compatible with the bornology induced by $E'$, and let 
$c:\Bbb R \to  E$ be a curve. Then the following conditions are 
equivalent.
\roster
\item The curve $c$ is topologically real analytic.
\item The curve $c$ is bornologically real analytic.
\item The curve $c$ extends to a holomorphic curve from some open 
neighborhood $U$ of $\Bbb R$ in $\Bbb C$ into the 
complexification $(E_{\Bbb C},E'_{\Bbb C})$.
\endroster
\endproclaim
\demo{Proof}
\therosteritem{1} $\Rightarrow$ \therosteritem{3}.
For every $t\in \Bbb R$ one has for some $\de>0$ and all $|s|<\de$ a 
converging power series representation 
$c(t+s)=\sum_{k=1}^{\infty}x_ks^k$. For any complex number $z$ with 
$|z|<\de$ the series converges in $(E_{\Bbb C},E'_{\Bbb C})$, hence $c$ 
can be locally extended to a holomorphic curve into $E_{\Bbb C}$. By the 
1-dimensional uniqueness theorem for holomorphic maps, these local 
extensions fit together to give a holomorphic extension as required.

\therosteritem{3} $\Rightarrow$ \therosteritem{2}.
A holomorphic curve factors locally over $(E_{\Bbb C})_B$,
where $B$ can be chosen of the form $B\x\sqrt{-1}B$. Hence the 
restriction of this factorization to $\Bbb R$ is real analytic into 
$E_B$.

\therosteritem{2} $\Rightarrow$ \therosteritem{1}.
Let $c$ be bornologically real analytic, i.e. $c$ is locally real 
analytic into some $E_B$, which we may assume to be complete. Hence $c$ 
is locally even topologically real analytic in $E_B$ by \nmb!{1.6} and 
so also in $E$.
\qed\enddemo

\proclaim{\nmb.{1.8}. Lemma} Let 
$E$ be a regular (i.e. every bounded set is contained and bounded in some 
step $E_\al$) inductive 
limit of complex locally convex spaces $E_\al\subseteq E$, let $c:\Bbb 
C\supseteq U\to  E$ be a holomorphic mapping, and let $W\subseteq \Bbb C$ be 
open and such that the closure $\overline W$ is compact and contained in $U$.
Then there exists some $\al$, such that $c|W:W \to  E_\al$ is well defined 
and holomorphic.
\endproclaim
\demo{Proof} Since $W$ is relatively compact, $c(W)$ is bounded in $E$.
It suffices to show that for the absolutely convex closed hull
$B$ of $c(W)$ the Taylor series of $c$ at each $z\in W$ 
converges in $E_B$, i.e. that $c|W:W \to  E_B$ is holomorphic.
This follows from the 

{\bf Vector valued Cauchy inequalities.} If $r>0$ is smaller
than the radius of convergence at $z$ of $c$ then
$$\frac{r^k}{k!}c^{(k)}(z) \in B$$
where $B$ is the closed absolutely convex
hull of $\{\,c(w): |w-z| = r\}$.
(By the Hahn-Banach theorem this follows directly from the scalar
valued case.)

Thus we get
$$\sum_{k=n}^m \left(\frac{w-z}r\right)^k\cdot \frac{r^k}{k!}c^{(k)}(z)
\in \sum_{k=n}^m \left(\frac{w-z}{r}\right)^k\cdot B$$
and so 
$$\sum_k \frac{c^{(k)}(z)}{k!}(w-z)^k$$ 
is convergent in $E_B$ for $|w-z|<r$.
Since $B$ is contained and bounded in 
some $E_\al$ one has $c|W:W \to  E_B=(E_\al)_B \to  E_\al$ is holomorphic.
\qed\enddemo

This proof also shows that holomorphic curves with values
in complex convenient vector spaces are topologically and bornologically 
holomorphic (compare with \nmb!{1.3}).

\proclaim{\nmb.{1.9}. Theorem (Linear real analytic mappings)}
Let $(E,E')$ be a complete dual pair. For 
any linear functional $\la:E \to  \Bbb R$ the following assertions are 
equivalent.
\roster 
\item $\la$ is bounded.
\item $\la\o c \in C^\om(\Bbb R,\Bbb R)$ for each real analytic $c:\Bbb
        R\to E$. 
\endroster
\endproclaim

This will be generalized in \nmb!{2.7} to non-linear mappings.

\demo{Proof} ($\Uparrow$). Let $\la$ satisfy
\therosteritem2 and suppose that 
there is a bounded sequence $(x_k)$ such that $\la(x_k)$ is
unbounded. By passing to a subsequence we may
suppose that $|\la(x_k)|>k^{2k}$. Let $a_k:= k^{-k}\,x_k$, then
$(r^k\,a_k)$ is bounded and  $(r^k\,\la(a_k))$ is unbounded for
all $r>0$.
Hence the curve $c(t) := \sum_{k=0}^\infty t^k\,a_k$ is 
given by a Mackey convergent power series. So $\la\o c$ is real analytic
and near $0$ we have $\la(c(t))=\sum_{k=0}^\infty b_k\,t^k$ for
some $b_k\in \Bbb R$.
But 
$$\la(c(t)) = \sum_{k=0}^N\la(a_k)t^k +
t^N\la\left(\sum_{k>N}a_kt^{k-N}\right)$$
and $t\mapsto
\sum_{k>N}a_kt^{k-N}$ is still a Mackey converging power series in $E$.
Comparing coefficients we see that $b_k = \la(a_k)$ and
consequently $\la(a_k)r^k$ is bounded for some $r>0$, a contradiction.

($\Downarrow$). Let $c:\Bbb R\to E$ be real analytic. By theorem \nmb!{1.5} 
the set 
$$\{\frac1{k!}\,c^{(k)}(a)\,r_k:a\in K, k\in\Bbb N\}$$ 
is bounded
for all compact sets $K\subset \Bbb R$ and for all sequences 
$(r_k)$ with $r_k\,t^k\to 0$ for all $t>0$. 
Since $c$ is smooth and bounded linear mappings are smooth 
(\cit!{5}, 2.4.4), the function $\la\o c$
is smooth and $(\la\o c)^{(k)}(a) = \la(c^{(k)}(a))$. By applying
\nmb!{1.5} we obtain that $\la\o c$ is real analytic.
\qed\enddemo

\proclaim{\nmb.{1.10}. Lemma}
Let $(E,E^1)$ and $(E,E^2)$ be two complete dual pairs with the same 
underlying vector space $E$. Then following statements are 
equivalent:
\roster
\item They have the same bounded sets.
\item They have the same smooth curves.
\item They have the same real analytic curves.
\endroster 
\endproclaim
\demo{Proof} \therosteritem 1 $\Leftrightarrow$ \therosteritem 2. 
This was shown in \cit!{14}.

\therosteritem 1 $\Rightarrow$ \therosteritem 3. This follows from 
\nmb!{1.5}, which shows that real analyticity is a bornological
concept.

\therosteritem 1 $\Leftarrow$ \therosteritem 3. This follows from 
\nmb!{1.9}.
\qed\enddemo

\proclaim{\nmb.{1.11}. Lemma} If a cone of linear maps $T_\al:(E,E')\to 
(E_\al,E_\al')$ between complete dual pairs generates the
bornology on $E$, then a curve $c:\Bbb R\to E$ is $C^\om$ resp. $C^\infty$ 
provided all the composites $T_\al\o c:\Bbb R\to E_\al$ are.
\endproclaim
\demo{Proof} The statement on the smooth curves is shown in 
\cit!{5}. That on the real analytic 
curves follows again from the bornological condition of \nmb!{1.5}.
\qed\enddemo

\heading\totoc \nmb0{2}. Real analytic mappings \endheading

Parts of \nmb!{2.1} to \nmb!{2.5} can be found in 
\cit!{2}.
For $x$ in any vector space $E$ let $x^k$ denote the element
$(x,\dots,x)\in E^k$.

\proclaim{\nmb.{2.1}. Lemma (Polarization formulas)} Let $f: E\x
\cdots\x E \to F$ be an 
$k$-linear symmetric mapping between vector spaces. Then
we have:
$$\align &f(\row x1k) = \tfrac 1{k!}\sum_{\row \ep1k = 0}^1
        (-1)^{k-\Si\ep_j}
        f\left(\left(x_0+\tsize\sum\ep_jx_j\right)^k\right).\tag1\\
&f(x^k)= \tfrac1{k!}\sum_{j=0}^k(-1)^{k-j}\tsize\binom kj f((a+jx)^k).\tag2\\
&f(x^k)= \tfrac{k^k}{k!}\sum_{j=0}^k(-1)^{k-j}\tsize\binom kj
        f((a+\tfrac jk x)^k).\tag3 \\
&f(x_1^0+\la x_1^1,\dots,x_k^0+\la x_k^1) = 
        \sum_{\row \ep1k = 0}^1 \la^{\Si\ep_j}
        f(x_1^{\ep_1},\dots,x_k^{\ep_k}).\tag4 
\endalign$$
\endproclaim
Formula \thetag4 will mainly be used for $\la=\sqrt{-1}$ in the
passage to the complexification.
\demo{Proof} \thetag1. (see \cit!{17}). By
multilinearity and symmetry the right hand side expands to
$$\sum_{j_0+\dots+j_k=k}\frac{A_{\row j0k}}{j_0!\cdots j_k!}
f(\undersetbrace j_0 \to{\row x00},\dots,\undersetbrace j_k\to
{\row xkk}),$$ where the coefficients are given by 
$$A_{\row j0k}=
\sum_{\row\ep1k=0}^1(-1)^{k-\Si\ep_j}\ep_1^{j_1}\cdots\ep_k^{j_k}.$$
The only nonzero coefficient is $A_{0,1,\dots,1}=1$.

\thetag2. In formula \thetag1 we put $x_0=a$ and all $x_j=x$.

\thetag3. In formula \thetag2 we replace $a$ by $ka$ and pull
$k$ out of the $k$-linear expression $f((ka+jx)^k)$.

\thetag4 is obvious.
\qed\enddemo
\proclaim{\nmb.{2.2}. Lemma (Power series)} Let $E$ be a real or
complex Fr\'echet 
space and let $f_k$ be a $k$-linear symmetric scalar valued
bounded functional on $E$, for each $k\in \Bbb N$. Then the
following statements are equivalent:
\roster
\item $\sum_kf_k(x^k)$ converges pointwise on an absorbing
        subset of $E$.
\item $\sum_kf_k(x^k)$ converges uniformly and absolutely on
        some neighborhood of 0.
\item $\{f_k(x^k): k\in\Bbb N, x\in U\}$ is bounded for some
        neighborhood $U$ of zero.
\item $\{f_k(\row x1k):k\in \Bbb N, x_j\in U\}$ is bounded for
        some neighborhood $U$ of 0.
\endroster
If any of these statements are satisfied over the reals, then also for the
complexification of the functionals $f_k$.
\endproclaim
\demo{Proof} \therosteritem1 $\Rightarrow$ \therosteritem3
The set $A_{K,r}:=\{x\in E:|f_k(x^k)|\leq Kr^k
\text{ for all }k\}$ is closed in $E$ since every bounded
multi linear mapping is continuous. The union $\bigcup_{K,r}A_{K,r}$ is 
$E$, since the series converges pointwise on an absorbing subset.
Since $E$ is Baire there are $K>0$ and $r>0$ such that the interior $U$ of
$A_{K,r}$ is non void. Let $x_0\in U$ and let $V$ be an
absolutely convex neighborhood of 0 contained in $U-x_0$

 From \nmb!{2.1} \thetag3 we get for all $x\in V$ the following estimate:
$$\align |f(x^k)| &\leq \tfrac{k^k}{k!}\sum_{j=0}^k\tsize\binom kj 
|f((x_0+\tfrac jk x)^k)| \\
&\leq \tfrac {k^k}{k!}2^kKr^k \leq K(2re)^k.\endalign$$
Now we replace $V$ by $\frac1{2re}\,V$ and get the result.

\therosteritem3 $\Rightarrow$ \therosteritem4.
 From \nmb!{2.1} \thetag1 we get for all $x_j\in U$ the estimate:
$$\align |f(\row x1k)| &\leq \tfrac 1{k!}\sum_{\row \ep1k = 0}^1
|f\left(\left(\tsize\sum\ep_jx_j\right)^k\right)| \\
&\leq \tfrac 1{k!}\sum_{\row \ep1k = 0}^1\left(\tsize\sum\ep_j\right)^k
|f\left(\left(\tfrac{\tsize\sum\ep_jx_j}{\tsize\sum\ep_j}\right)^k\right)|\\
&\leq \tfrac 1{k!}\sum_{\row \ep1k = 0}^1
        \left(\tsize\sum\ep_j\right)^k C\\
&\leq \tfrac 1{k!}\sum_{j=0}^k\tsize\binom kj j^k \,C
        \leq C(2e)^k.\endalign$$
Now we replace $U$ by $\frac1{2e}\,U$ and get \therosteritem4.

\therosteritem4 $\Rightarrow$ \therosteritem2. The series
converges on $rU$ uniformly and absolutely for any $0<r<1$.

\therosteritem2 $\Rightarrow$ \therosteritem1 is clear.

\therosteritem4, real case, $\Rightarrow$ \therosteritem4,
complex case, by \nmb!{2.1}.\thetag4 for $\la=\sqrt{-1}$.
\qed\enddemo

\proclaim{\nmb.{2.3}. Theorem (Holomorphic functions on
Fr\'echet spaces)} Let $U\subseteq E$ be open in a
complex Fr\'echet space $E$. The following statements on 
$f:U\to \Bbb C$ are equivalent:
\roster
\item $f$ is holomorphic along holomorphic curves.
\item $f$ is smooth and the derivative $df(z):E\to \Bbb C$ 
        is $\Bbb C$-linear for all $z\in U$.
\item $f$ is smooth and is locally given by its pointwise
        converging Taylor series.
\item $f$ is smooth and is locally given by its uniformly and absolutely
        converging Taylor series.
\item $f$ is locally given by a uniformly and absolutely
        converging power series.
\endroster
\endproclaim
\demo{Proof} \therosteritem1 $\Leftrightarrow$
\therosteritem2 \cit!{15}, 2.12.

\therosteritem1 $\Rightarrow$ \therosteritem3.
Let $z\in U$ be arbitrary, without loss of generality 
$z=0$, and let $b_n:=\frac{f^{(n)}(z)}{n!}$ be the n-th Taylor
coefficient of $f$ at $z$. Then $b_n:E^n\to \Bbb C$
is symmetric, n-linear and bounded and the series
$\sum_{n=0}^{\infty}b_n(v,..,v)t^n$ converges to $f(z+tv)$ for
small $t$. Hence the set of those $v$ for which the series
$\sum_{n=0}^{\infty}b_n(v,..,v)$ converges is absorbing. 
By \nmb!{2.2}, \therosteritem1 $\Rightarrow$ \therosteritem2 it
converges on a neighborhood of 0 to $f(z+v)$.

\therosteritem3 $\Rightarrow$ \therosteritem4 follows from
\nmb!{2.2},\therosteritem2 $\Rightarrow$ \therosteritem3.

\therosteritem4 $\Rightarrow$ \therosteritem5 is obvious.

\therosteritem5 $\Rightarrow$ \therosteritem1 is the chain
rule for converging power series, which easily can be shown using  
\nmb!{2.2}, \therosteritem2 $\Rightarrow$ \therosteritem4.
\qed\enddemo

\proclaim{\nmb.{2.4}. Theorem (Real analytic functions on
Fr\'echet spaces)} Let $U\subseteq E$ be open in a
real Fr\'echet space $E$. The following statements on 
$f:U\to \Bbb R$ are equivalent:
\roster
\item $f$ is smooth and is real analytic along topologically
        real analytic curves. 
\item $f$ is smooth and is real analytic along affine lines.
\item $f$ is smooth and is locally given by its pointwise
        converging Taylor series.
\item $f$ is smooth and is locally given by its uniformly and absolutely
        converging Taylor series.
\item $f$ is locally given by a uniformly and absolutely
        converging power series.
\item $f$ extends to a holomorphic mapping $\tilde f:\tilde U\to \Bbb C$
        for an open subset $\tilde U$ in the complexification
        $E_{\Bbb C}$ with $\tilde U\cap E=U$.
\endroster
\endproclaim

\demo{Proof} \therosteritem1 $\Rightarrow$ \therosteritem2
is obvious.

\therosteritem2 $\Rightarrow$ \therosteritem3.
Repeat the proof of \nmb!{2.3}, \therosteritem1 $\Rightarrow$
\therosteritem3.

\therosteritem3 $\Rightarrow$ \therosteritem4 follows from
\nmb!{2.2},\therosteritem2 $\Rightarrow$ \therosteritem3.

\therosteritem4 $\Rightarrow$ \therosteritem5 is obvious.

\therosteritem5 $\Rightarrow$ \therosteritem6. Locally we
can extend converging power series into the complexification by
\nmb!{2.2}. Then we take the union $\tilde U$ of their domains of definition
and use uniqueness to glue $\tilde f$ which is holomorphic by
\nmb!{2.3}.

\therosteritem6 $\Rightarrow$ \therosteritem1. Obviously $f$
is smooth. Any topologically real analytic curve $c$ in $E$ can locally be extended to
a holomorphic curve in $E_{\Bbb C}$ by \nmb!{1.3}. So $f\o c$ is
real analytic.
\qed\enddemo

\subheading{\nmb.{2.5}}
The assumptions "$f$ is smooth" cannot be dropped in \nmb!{2.4}.1 
even in finite dimensions, as shown by the following example,
due to \cit!{3}.

\proclaim{Example} The mapping $f:\Bbb R^2\to \Bbb R$, defined by
$$f(x,y):=\frac{xy^{n+2}}{x^2+y^2}$$
is real analytic along real analytic
curves, is n-times continuous differentiable but is not smooth and hence not
real analytic.
\endproclaim
\demo{Proof}
Take a real analytic curve $t\mapsto (x(t),y(t))$ into $\Bbb R^2$. The
components can be factored as $x(t)=t^nu(t)$, $y(t)=t^nv(t)$
for some $n$
and real analytic curves $u$, $v$ with $u(0)^2+v(0)^2\not= 0$. The
composite $f\o (x,y)$ is
then the function 
$$t\mapsto t^n\frac{uv^{n+2}}{u^2+v^2}(t),$$
which is obviously real analytic near $0$.
The mapping $f$ is n-times continuous differentiable, since it is
real analytic on $\Bbb R^2\backslash\{0\}$ and the directional
derivatives of order $i$ are $(n+1-i)$-homogeneous, hence
continuously extendable to $\Bbb R^2$.
But $f$ cannot be $(n+1)$-times continuous differentiable, 
otherwise the derivative of order $n+1$ would be 
constant, and hence $f$ would be a polynomial.
\qed\enddemo

\subheading{\nmb.{2.6}. Definition (Real analytic mappings)} Let
$(E,E')$ be a dual pair. 
Let us denote by $C^\om(\Bbb R,E)$ the \idx{\it space of all real analytic 
curves}.

Let $U\subseteq E$ be $c^\infty$-open, and let $(F,F')$ be a second dual 
pair.
A mapping $f:U\to F$ will be called \idx{\it real analytic} or
$C^\om$ for short, if $f$ is real analytic along real analytic
curves and is smooth  
(i.e. is smooth along smooth curves); so $f\o c \in C^\om(\Bbb R,F)$
for all $c\in C^\om(\Bbb R,E)$ with $c(\Bbb R)\subseteq U$
and $f\o c \in C^\infty(\Bbb R,F)$ for all $c\in C^\infty(\Bbb R,E)$ with
$c(\Bbb R)\subseteq U$.
Let us denote by $C^\om(U,F)$ the \idx{\it space of all real analytic  
mappings} from $U$ to $F$.

\proclaim{\nmb.{2.7}. Hartogs' Theorem for real analytic
mappings} Let $(E,E')$ and $(F,F')$ be 
complete dual pairs, let $U\subseteq E$ be $c^\infty$-open, and
let $f: U\to F$.
Then $f$  is real analytic if and only if $f$ is smooth and 
$\la\o f$ is real analytic along each affine line in $E$, for
all $\la\in F'$.
\endproclaim
\demo{Proof} One direction is clear, and by definition \nmb!{2.6}
we may assume that $F=\Bbb R$. 

Let $c:\Bbb R\to U$ be real analytic. 
We show that $f\o c$ is real analytic by using lemma 
\nmb!{1.5}. So let $(r_k)$ be a sequence such that
$r_kr_\ell\geq r_{k+\ell}$ and $r_k\,t^k\to 0$ for all $t>0$ and let 
$K\subset \Bbb R$ be compact.
We have to show, that there is an $\ep>0$ such that the set
$$\left\{\tfrac1{\ell!}(f\o c)^{(\ell)}(a)\,r_l\,(\tfrac\ep2)^\ell:a\in K, 
\ell\in \Bbb N\right\}$$ 
is bounded.

By theorem \nmb!{1.5} the set 
$$\left\{\tfrac1{n!}c^{(n)}(a)\,r_n:n\geq1, a\in K\right\}$$
is contained in some boun\-ded absolutely convex subset
$B\subseteq E$, such that $E_B$ 
is a Banach space.
Clearly for the inclusion $i_B:E_B\to E$ the function $f\o i_B$
is smooth and real analytic along affine lines. Since $E_B$ is a
Banach space, by \nmb!{2.4}, \therosteritem2 $\Rightarrow$
\therosteritem4 $f\o i_B$ is locally given by its uniformly and
absolutely converging Taylor series. Then by \nmb!{2.2},
\therosteritem2 $\Rightarrow$ \therosteritem4
there is an $\ep>0$ such that the set
$$\left\{\tfrac1{k!}d^kf(c(a))(\row x1k): k\in \Bbb N, 
x_j\in \ep B,a\in K\right\}$$
is bounded, so is contained in $[-C,C]$ for some $C>0$.

The Taylor series of $f\o c$ at $a$
is given by 
$$(f\o c)(a+t) = \sum_{\ell\geq0} 
\sum_{k\geq0}\tfrac1{k!}
\sum\Sb (m_n)\in \Bbb N_0^{\Bbb N} \\
        \sum_n m_n = k \\
        \sum_n m_n\,n = \ell \endSb
        \frac {k!}{\prod_n m_n!}
        d^kf(c(a))\Bigl( \prod_n (\tfrac1{n!}c^{(n)}(a))^{m_n}
\Bigr) t^\ell,
$$
where
$$\prod_n x_n^{m_n} := 
        (\undersetbrace m_1 \to{\row x11},\dots,
        \undersetbrace m_n\to {\row xnn},\dotsc).$$
This follows easily from composing the Taylor series of $f$ and
$c$ and ordering by powers of $t$. Furthermore we have
$$\sum\Sb (m_n)\in \Bbb N_0^{\Bbb N} \\
        \sum_n m_n = k \\
        \sum_n m_n\,n = \ell \endSb
    \frac {k!}{\prod_n m_n!} = \tsize\binom{\ell-1}{k-1}$$
by the following argument: It is the
$\ell$-th Taylor coefficient at 0 of the function
$$\left(\sum_{n\geq0}t^n-1\right)^k = \left(\frac t{1-t}\right)^k =
t^k\sum_{j=0}^\infty \binom {-k}j(-t)^j,$$
which turns out to be
the binomial coefficient in question.

By the foregoing considerations we may estimate as follows.
$$\allowdisplaybreaks\align
&\tfrac1{\ell!}|(f\o c)^{(\ell)}(a)|\,r_l\,(\tfrac\ep2)^\ell\leq \\
&\leq \sum_{k\geq0}\Bigl|\tfrac1{k!}
\sum\Sb (m_n)\in \Bbb N_0^{\Bbb N} \\
        \sum_n m_n = k \\
        \sum_n m_n\,n = \ell \endSb
        \frac {k!}{\prod_n m_n!} 
d^kf(c(a))\Bigl(\prod_n (\tfrac1{n!}c^{(n)}(a))^{m_n}
\Bigr)\Bigr|\,r_\ell\, (\tfrac\ep2)^\ell \\
&\leq  \sum_{k\geq0}\Bigl|\tfrac1{k!}
\sum\Sb (m_n)\in \Bbb N_0^{\Bbb N} \\
        \sum_n m_n = k \\
        \sum_n m_n\,n = \ell \endSb
    \frac {k!}{\prod_n m_n!} 
d^kf(c(a))\Bigl(\prod_n (\tfrac1{n!}c^{(n)}(a)\,r_n\,\ep^n)^{m_n}
\Bigr)\Bigr|\tfrac1{2^\ell} \\
&\leq \sum_{k\geq0}\tsize\binom {\ell-1}{k-1}\,C\tfrac1{2^\ell}
= \tfrac12 C,\endalign$$
because
$$\sum\Sb (m_n)\in \Bbb N_0^{\Bbb N} \\
        \sum_n m_n = k \\
        \sum_n m_n\,n = \ell \endSb
    \frac {k!}{\prod_n m_n!}\prod_n(\tfrac1{n!}c^{(n)}(a)\,\ep^n\,r_n)^{m_n}\,
\in \tsize\binom{\ell-1}{k-1}\,(\ep B)^k\subseteq(E_B)^k.\qed$$
\enddemo

\proclaim{\nmb.{2.8}. Corollary} Let $(E,E')$ and $(F,F')$ be
complete dual pairs, let $U\subseteq E$ be $c^\infty$-open, and
let $f: U\to F$.
Then $f$  is real analytic if and only if $f$ is smooth and 
$\la\o f\o c$ is real analytic for every periodic (topologically) real 
analytic curve $c:\Bbb R\to
U\subseteq E$ and all $\la\in F'$.
\endproclaim
\demo{Proof} By \nmb!{2.7} $f$ is real analytic  if and only if $f$ is smooth and $\la \o 
f$ is real analytic along topologically real analytic curves $c:\Bbb R\to  
E$. 
Let $h:\Bbb R \to  \Bbb R$ be defined by $h(t)=t_0+\ep\cdot \sin t$. Then 
$c\o h:\Bbb R \to  \Bbb R \to  U$ is 
a (topologically) real analytic, periodic function with period $2\pi$, 
provided $c$ is (topologically) real analytic. If $c(t_0)\in U$ we 
can choose $\ep > 0$ such that $h(\Bbb R)\subseteq c^{-1}(U)$. Since 
{\sl sin} is locally around 0 invertible, real analyticity of $\la\o f\o c 
\o h$ implies that $\la\o f\o c$ is real analytic near $t_0$. Hence the 
proof is completed.
\qed\enddemo

\proclaim{\nmb.{2.9}. Corollary (Reduction to Banach spaces)}
Let $(E,E')$ be a complete 
dual pair, let $U\subset E$ be $c^\infty$-open, and let $f:U\to
\Bbb R$ be a mapping. Then $f$ is real analytic if and only if
the restriction $f: E_B\supset U\cap E_B\to \Bbb R$ is real
analytic for all bounded absolutely convex subsets $B$ of $E$. 
\endproclaim

So any result valid on Banach spaces can be translated into a
result valid on complete dual pairs.

\demo{Proof} By theorem \nmb!{2.7} it suffices to check $f$ along 
bornologically real analytic curves. These factor by definition locally
to real analytic curves into some $E_B$.
\qed\enddemo

\proclaim{\nmb.{2.10}. Corollary} Let $U$ be a $c^\infty$-open
subset in a complete dual pair $(E,E')$ and let $f:U\to \Bbb R$
be real analytic. Then for every bounded $B$ there is some
$r_B>0$ such that the Taylor series 
$$y\mapsto \sum\tfrac1{k!}d^kf(x)(y^k)$$ 
converges to $f(x+y)$ uniformly and absolutely on $r_BB$.
\endproclaim
\demo{Proof} Use \nmb!{2.9} and \nmb!{2.4}.\therosteritem4.
\qed\enddemo

\subheading{\nmb.{2.11}} Scalar analytic functions on 
convenient vector spaces $E$ are in general not germs of holomorphic 
functions from $E_{\Bbb C}$ to $\Bbb C$:

\proclaim{Example} Let $f_k:\Bbb R \to  \Bbb R$ be real analytic functions 
with radius of convergence at zero converging to 0 for $k \to  \infty$.
Let $f:\Bbb R^{(\Bbb N)} \to  \Bbb R$ be the mapping defined on the 
countable sum $\Bbb R^{(\Bbb N)}$ of the reals by 
$f(x_0,x_1,...):=\sum_{k=1}^{\infty}x_kf_k(x_0)$. Then $f$ is real 
analytic, but there is no complex valued holomorphic mapping $\tilde f$ 
on some neighborhood of 0 in $\Bbb C^{(\Bbb N)}$ which extends
$f$, and the Taylor series of $f$ is not pointwise convergent on any
$c^\infty$-open neighborhood of 0.
\endproclaim

\demo{Proof}
{\bf Claim.} $f$ is real analytic. \newline
Since the limit $\Bbb R^{(\Bbb N)} = \varinjlim_n \Bbb R^n$ is regular, 
every smooth curve (and hence every real analytic curve) in 
$\Bbb R^{(\Bbb N)}$ is locally smooth (resp. real analytic) 
into $\Bbb R^n$ for some $n$. Hence $f\o c$ is locally just a finite 
sum of smooth (resp. real analytic) functions and is therefore smooth 
(resp. real analytic).

{\bf Claim.} $f$ has no holomorphic extension.\newline 
Suppose there exists some 
holomorphic extension $\tilde f:U\to \Bbb C$, where $U\subseteq 
\Bbb C^{(\Bbb N)}$ is $c^\infty$-open neighborhood of 0, and is therefore 
open in the locally convex Silva topology by 
\cit!{5}, 6.1.4.ii. Then $U$ is even open
in the box-topology  
\cit!{10}, 4.1.4, i.e\. there exist $\ep_k>0$ for all $k$, 
such that $\{(z_k)\in \Bbb C^{(\Bbb N)}: 
|z_k| \leq \ep_k \text{ for all }k \} \subseteq U$. Let $U_0$ be the open 
disk in $\Bbb C$ with radius $\ep_0$ and let $\tilde f_k: U_0 \to  \Bbb C$
be defined by 
$\tilde f_k(z):=\tilde f(z,0,...,0,\ep_k,0,...)\tfrac1{\ep_k}$, where 
$\ep_k$ is inserted instead of the variable $x_k$. Obviously $\tilde 
f_k$ is an extension of $f_k$, which is impossible, since the radius of 
convergence of $f_k$ is less than $\ep_0$ for $k$ sufficiently large.

{\bf Claim.} The Taylor series does not converge. \newline
If the Taylor series would be pointwise convergent on some $U$,
then the previous arguments would show that the radii of
convergence of the $f_k$ were bounded from below.
\qed\enddemo

\heading\totoc \nmb0{3}. Function spaces in finite dimensions \endheading

\subheading{\nmb.{3.1}. Spaces of holomorphic functions}
For a complex manifold $N$ (always assumed to be separable) let
$\Cal H(N,\Bbb C)$ be the \idx{\it space 
of all holomorphic} functions on $N$ with the topology
of uniform convergence on compact subsets of $N$.

Let $\Cal H_b(N,\Bbb C)$ denote the Banach space of
bounded holomorphic functions on $N$ equipped with the supremum norm. 

For any open subset $W$ of $N$ let $\Cal H_{bc}(W\subseteq
N,\Bbb C)$ be the closed subspace of $\Cal H_b(W,\Bbb C)$ of all
holomorphic functions on $W$ which extend to continuous
functions on the closure $\overline W$. 

For a poly-radius  $r=(\row r1n)$ with $r_i>0$ and for $1\le p\le\infty$ 
let $\ell^p_r$ denote the real Banach space 
$\left\{\;x\in\Bbb R^{\Bbb N^n}:
\|(x_\al r^\al)_{\al \in \Bbb N^n}\|_p < \infty\;\right\}$.

\proclaim{\nmb.{3.2}. Theorem (Structure of $\Cal H(N,\Bbb C)$ for 
complex manifolds $N$)}\newline
The space $\Cal H(N,\Bbb C)$ 
of all holomorphic functions on $N$ with the topology
of uniform convergence on compact subsets of $N$ is a (strongly) nuclear 
Fr\'echet space and embeds as a closed subspace into 
$C^\infty(N,\Bbb R)^2$.
\endproclaim 
\demo{Proof}
By taking a countable covering of $N$ with compact sets, one
obtains a countable  
neighborhood basis of 0 in $\Cal H(N,\Bbb C)$. Hence $\Cal H(N,\Bbb C)$ 
is metrizable.

That $\Cal H(N,\Bbb C)$ is complete, and hence a Fr\'echet space, 
follows since the limit of a sequence of holomorphic functions with 
respect to the topology of uniform convergence on compact sets is again 
holomorphic.

The vector space $\Cal H(N,\Bbb C)$ is a subspace of 
$C^\infty(N,\Bbb R^2)=C^\infty(N,\Bbb R)^2$ since a function 
$N\to\Bbb C$ is holomorphic if and only if it is smooth and the derivative at every 
point is $\Bbb C$-linear. It is a closed subspace, since it is
described by the continuous linear equations
$df(x)(\sqrt{-1}\cdot v)=\sqrt{-1}\cdot df(x)(v)$.
Obviously the identity from $\Cal H(N,\Bbb C)$ with the subspace 
topology to $\Cal H(N,\Bbb C)$ is continuous, hence by the open mapping 
theorem \cit!{10}, 5.5.2, for Fr\'echet spaces it is an isomorphism.

That $\Cal H(N,\Bbb C)$ is nuclear and unlike $C^\infty(N,\Bbb R)$ even 
strongly nuclear can be shown as follows. For $N$ equal to the open 
unit disk $\Bbb D\subseteq \Bbb C$ this result can be found in 
\cit!{10}, 21.8.3.b. More generally for $N=\Bbb D^n$ one has 
that $\Cal H(N,\Bbb C)=\bigcap_{0<r<1}\ell^1_{(r,..,r)}\otimes\Bbb C$ as 
vector spaces. The identity from the right to the left is obviously 
continuous, if the intersection is supplied with the projective limit 
topology induced from the Banach spaces $\ell^1_{(r,..,r)}\otimes\Bbb 
C$, a Fr\'echet topology. Hence again by the open mapping theorem it is 
an isomorphism.
Using now the Grothendieck-Pietsch criterion, cf. \cit!{10}, 21.8.2, 
one concludes that $\Cal H(\Bbb D^n,\Bbb C)$ is strongly nuclear,
see also \cit!{30}, p. 530.
For an arbitrary $N$ the space $\Cal H(N,\Bbb C)$ carries the 
initial topology induced by the linear mappings 
$u_*:\Cal H(N,\Bbb C) \to \Cal H(u(U),\Bbb C)$ for all charts $(u,U)$ of 
$N$, for which we may assume $u(U)=\Bbb D^n$, and hence by the 
stability properties of strongly nuclear spaces, cf. \cit!{10}, 21.1.7, 
$\Cal H(N,\Bbb C)$ is strongly nuclear.
\qed\enddemo

\subheading{\nmb.{3.3} Spaces of germs of holomorphic functions}
For a subset $A\subseteq N$ let 
$\Cal H(A\subseteq N, \Bbb C)$ 
be the space of \idx{\it germs along $A$ of holomorphic functions}  
$W\to \Bbb C$ for open sets $W$ in $N$ containing $A$. 
We equip $\Cal H(A\subseteq N,\Bbb C)$ with the locally convex topology
induced by the inductive cone 
$\Cal H(W,\Bbb C)\to \Cal H(A\subseteq N,\Bbb C)$ for all $W$.
This is Hausdorff, since iterated derivatives at points in $A$
are continuous functionals and separate points.
In particular $\Cal H(W\subseteq N,\Bbb C)=\Cal H(W,\Bbb C)$ for
$W$ open in $N$.
For $A_1\subset A_2\subset N$ the "restriction" mappings
$\Cal H(A_2\subset N,\Bbb C)\to \Cal H(A_1\subset N,\Bbb C)$ are
continuous.

The structure of $\Cal H(A\subseteq S^2,\Bbb C)$,
where $A\subseteq S^2$ is a subset of the Riemannian sphere, has
been studied by \cit!{29}, 
\cit!{26}, \cit!{31}, 
\cit!{12}, and \cit!{8}.

\proclaim{\nmb.{3.4}. Theorem (Structure of $\Cal H(K\subseteq N,\Bbb C)$ 
for compact subsets $K$ of complex manifolds $N$)}
The following inductive cones are cofinal to each other.
$$\gather
\Cal H(K\subseteq N,\Bbb C) 
    \gets \{\Cal H(W,\Bbb C),N\supseteq W\supseteq K\}\\
\Cal H(K\subseteq N,\Bbb C) 
    \gets \{\Cal H_b(W,\Bbb C),N\supseteq W\supseteq K\} \\
\Cal H(K\subseteq N,\Bbb C) \gets 
    \{\Cal H_{bc}(W\subseteq N,\Bbb C),N\supseteq W\supseteq K\}\\ 
\endgather$$

If $K = \{z\}$ these inductive cones and the following ones for
$1\le p\le \infty$ are cofinal to each other.
$$\Cal H(\{z\}\subseteq N,\Bbb C) \gets 
\{\ell^p_r\otimes\Bbb C,\;r\in \Bbb R_+^n\}$$

So all inductive limit topologies coincide.
Furthermore, the space $\Cal H(K\subseteq N,\Bbb C)$ is a Silva space, 
i.e. a countable inductive limit of Banach spaces, where the connecting 
mappings between the steps are compact, i.e. mapping bounded sets to 
relatively compact ones. The connecting mappings are even strongly 
nuclear. In particular, the limit is regular, i.e. every 
bounded subset is contained and bounded in some step, and $\Cal 
H(K\subseteq N,\Bbb C)$ is complete and (ultra-)bornological (hence a 
convenient vector space), webbed, strongly nuclear, 
reflexive and its dual is a strongly nuclear Fr\'echet space. It is 
however not a Baire space.
\endproclaim
\demo{Proof} Let $K\subseteq V \subseteq \overline V \subseteq W
\subseteq N$, where
$W$ and $V$ are open and $\overline V$ is compact. Then the obvious
mappings 
$$\Cal H_{bc}(W\subseteq N,\Bbb C)\to \Cal H_{b}(W,\Bbb C)\to 
\Cal H(W,\Bbb C) \to  \Cal H_{bc}(V\subseteq N,\Bbb C)$$
are continuous. This implies the first cofinality assertion. For $q\le p$
and $s<r$ the obvious maps $\ell^q_r \to  \ell^p_r$,
$\ell^\infty_r \to  \ell^1_s$, and $\ell^1_r\otimes\Bbb C \to  \Cal
H_{b}(\{w\in \Bbb C^n:|w_i-z_i|<r_i\},\Bbb C)\to  \ell^\infty_s\otimes 
\Bbb C$ 
are continuous, by
the Cauchy inequalities. So the remaining cofinality assertion follows.

Let us show next that the connecting mapping
$\Cal H_{b}(W,\Bbb C) \to  \Cal H_{b}(V,\Bbb C)$ is strongly nuclear (hence 
nuclear and compact). 
Since the restriction mapping from $E:=\Cal H(W,\Bbb C)$ to 
$\Cal H_b(V,\Bbb C)$ is continuous, it factors over 
$E \to  \widetilde {E_{(U)}}$ for some 0-neigh\-bor\-hood $U$ in $E$, where
$\widetilde {E_{(U)}}$ is the completed quotient of $E$ with the
Minkowski functional of $U$ as norm, see \cit!{10}, 6.8.
Since $E$ is strongly nuclear by \nmb!{3.2}, there exists by definition 
some larger 0-neighborhood $U'$ in $E$ such that the natural mapping 
$\widetilde {E_{(U')}} \to  \widetilde {E_{(U)}}$ is strongly nuclear. So the 
claimed connecting mapping is strongly nuclear, since it can be 
factorized as
$$\Cal H_{b}(W,\Bbb C) \to  
\Cal H(W,\Bbb C)=E \to 
\widetilde {E_{(U')}} \to  
\widetilde {E_{(U)}} \to 
\Cal H_{b}(V,\Bbb C)
$$

That a Silva space is regular and complete, can be found in 
\cit!{4}, 7.4 and 7.5.

That $\Cal H(K\subseteq N,\Bbb C)$ is ultra-bornological, webbed and 
strongly nuclear follows from
the permanence properties of ultra-bornological spaces, 
\cit!{10}, 13.2.5, of 
webbed spaces \cit!{10}, 5.3.3 
and of strongly nuclear spaces \cit!{10}, 21.1.7.

Furthermore, $\Cal H(K\subseteq N,\Bbb C)$ is reflexive and its strong 
dual is a Fr\'echet space, since it is a Silva-space, cf. 
\cit!{10}, 12.5.9 and p.270. The dual is even strongly nuclear, since 
$\Cal H(K\subseteq N,\Bbb C)$ is a nuclear Silva-space, cf. 
\cit!{10}, 21.8.6. 

The space $\Cal H(K\subseteq N,\Bbb C)$ has however 
not the Baire property, since it is webbed but not metrizable, cf. 
\cit!{10}, 5.4.4. If it were metrizable then it would
be of finite dimension, by \cit!{4}, 7.7. This is not
the case.
\qed\enddemo

Completeness of $\Cal H(K\subseteq \Bbb C^n,\Bbb C)$ was shown in 
\cit!{31}, th\'eor\`eme II, 
and regularity of the inductive limit 
$\Cal H(K\subseteq \Bbb C,\Bbb C)$ can be found in 
\cit!{12}, Satz 12.

\proclaim{\nmb.{3.5}. Lemma}
For a closed subset $A\subseteq \Bbb C$ the spaces $\Cal H(A\subseteq 
S^2,\Bbb C)$
and the space
$\Cal H_\infty(S^2\setminus A\subseteq S^2,\Bbb C)$ of all germs
vanishing at $\infty$ are strongly dual to each
other.
\endproclaim
\demo{Proof}
This is due to \cit!{12}, Satz 12 and has been generalized 
by \cit!{8}, th\'eor\`eme 2 bis, 
to arbitrary subsets $A\subseteq S^2$. 
\qed\enddemo 

Compare also the modern theory of 
hyperfunctions, cf. \cit!{11}.

\proclaim{\nmb.{3.6}. Theorem (Structure of $\Cal H(A\subseteq N,\Bbb C)$ 
for closed subsets $A$ of complex manifolds $N$)}
The inductive cone
$$\Cal H(A\subseteq N,\Bbb C) \gets 
\{\;\Cal H(W,\Bbb C): A\subseteq W \underset\text{open}\to\subseteq N\}$$
is regular, i.e. every bounded set is contained and bounded in
some step.

The projective cone
$$\Cal H(A\subseteq N,\Bbb C) \to \{\;\Cal H(K\subseteq N,\Bbb C):
K\text{ compact in } A \}$$
generates the bornology of $\Cal H(A\subseteq N,\Bbb C)$.

The space $\Cal H(A\subseteq N,\Bbb C)$ is Montel (hence
quasi-complete and reflexive), and ultra-bornological (hence a convenient
vector space). Furthermore it is webbed and conuclear.
\endproclaim
\demo{Proof} Compare also with the proof of the more general theorem 
\nmb!{7.3}.

We choose a continuous function $f:N\to\Bbb R$
which is positive and proper.
Then $(f\i([n,n+1]))_{n\in\Bbb N_0}$ is an exhaustion of $N$ by
compact subsets and 
$$(K_n:=A\cap f\i([n,n+1]))$$ 
is a compact exhaustion of $A$.

Let $\Cal B\subseteq \Cal H(A\subseteq N,\Bbb C)$ be bounded. Then
$\Cal B|K$ is also bounded in 
$\Cal H(K\subseteq N,\Bbb C)$ for each compact subset $K$ of $A$.
Since the cone
$$\{\Cal H(W,\Bbb C):K\subseteq W \underset\text{open}\to\subseteq N\}\to
\Cal H(K\subseteq N,\Bbb C)$$
is regular by \nmb!{3.4}, there exist open subsets $W_K$ of $N$
containing $K$ such that $\Cal B|K$ is contained (so that the
extension of each germ is unique) and bounded in
$\Cal H(W_K,\Bbb C)$. In particular we choose 
$W_{K_n\cap K_{n+1}}\subseteq W_{K_n}\cap W_{K_{n+1}}\cap f\i((n,n+2))$.
Then we put
$$W:=\bigcup_n(W_{K_n}\cap f\i((n,n+1)))\cup \bigcup_nW_{K_n\cap K_{n+1}}.$$
It is easily checked that $W$ is open in $N$, contains $A$, and
that each germ in $\Cal B$ has a unique extension to $W$.
$\Cal B$ is bounded in $\Cal H(W,\Bbb C)$ if it is uniformly
bounded on each compact subset $K$ of $W$. Each $K$ is
covered by finitely many $W_{K_n}$ and $\Cal B|K_n$ is bounded
in $\Cal H(W_{K_n},\Bbb C)$, so $\Cal B$ is bounded as required.

The space $\Cal H(A\subseteq N,\Bbb C)$ is ultra-bornological, Montel and 
in particular quasi-complete, and conuclear, as regular inductive 
limit of the nuclear Fr\'echet spaces $\Cal H(W,\Bbb C)$.

And it is webbed because it is the (ultra-)bornologification of the 
countable projective limit of webbed spaces 
$\Cal H(K\subseteq N,\Bbb C)$, cf. \cit!{10}, 13.3.3 + 5.3.3. 
\qed\enddemo

\proclaim{\nmb.{3.7}. Lemma} Let $A$ be closed in $\Bbb C$. Then the dual
generated by the projective cone
$$\Cal H(A\subseteq\Bbb C,\Bbb C) \to 
\{\;\Cal H(K\subset\Bbb C,\Bbb C),\;K\text{ compact in }A\;\}
$$
is just the topological dual of $\Cal H(A\subseteq\Bbb C,\Bbb C)$. 
\endproclaim
\demo{Proof} The induced topology is obviously coarser than the
given one. So let $\la$ be a continuous linear functional on
$\Cal H(A\subseteq\Bbb C,\Bbb C)$. 
Then we have 
$\la \in \Cal H_\infty(S^2\setminus A\subseteq S^2,\Bbb C)$ by \nmb!{3.5}.
Hence $\la \in \Cal H(U,\Bbb C)$ for some open neighborhood $U$
of $S^2\setminus A$, so again by \nmb!{3.5} $\la$ is a continuous
functional on $\Cal H(K\subset S^2,\Bbb C)$, where $K = S^2\setminus U$ is
compact in $A$. So $\la$ is continuous for the induced topology.
\qed\enddemo

\proclaim{Problem} Does this cone generate even the
topology of $\Cal H(A\subseteq \Bbb C,\Bbb C)$?
This would imply that the bornological topology on $\Cal
H(A\subseteq\Bbb C,\Bbb C)$ is complete and nuclear.
\endproclaim

\proclaim{\nmb.{3.8}. Lemma (Structure of $\Cal H(A\subseteq N,\Bbb C)$ 
for smooth closed submanifolds $A$ of complex manifolds $N$)}
The projective cone
$$\Cal H(A\subseteq N,\Bbb C) \to \{\;\Cal H(\{z\}\subseteq N,\Bbb C):
z\in A\}$$
generates the bornology.
\endproclaim
\demo{Proof}
Let $\Cal B\subseteq \Cal H(A\subseteq N,\Bbb C)$ be such that
the set $\Cal B$ is bounded in $\Cal H(\{z\}\subseteq N,\Bbb C)$
for all $z\in A$. By the  
regularity of the inductive cone 
$\Cal H(\{0\}\subseteq\Bbb C^n,\Bbb C) \leftarrow
\Cal H(W,\Bbb C)$ we find arbitrary small open neighborhoods $W_z$ such 
that the set $\Cal B_z$ of the germs at $z$ of all germs in $\Cal B$ is 
contained and bounded in $\Cal H(W_z,\Bbb C)$.

Now choose a tubular neighborhood $p:U\to A$ of $A$ in $N$. We may assume 
that $W_z$ is contained in $U$, has fibers which are star shaped with 
respect to the zero-section and the intersection with A is connected. 
The union $W$ of all the $W_z$, is therefore an open subset of $U$ 
containing $A$. And it remains to show that the germs in $\Cal B$ extend 
to $W$. For this it is enough to show that the extensions of the germs 
at $z_1$ and $z_2$ agree on the intersection of $W_{z_1}$ with 
$W_{z_2}$. So let $w$ be a point in the intersection. It can be radially 
connected with the base point 
$p(w)$, which itself can be connected by curves in A with $z_1$ and 
$z_2$. Hence the extension of both germs to $p(w)$ coincide with the 
original germ, and hence their extensions to $w$ are equal.

That $\Cal B$ is bounded in $\Cal H(W,\Bbb C)$, follows immediately 
since every compact subset $K\subseteq W$ can be covered by finitely many 
$W_z$.
\qed\enddemo

\subheading{\nmb.{3.9}}
The following example shows that \nmb!{3.8} fails to be true for 
general closed subsets $A\subseteq N$.

\proclaim{Example}
Let $A:=\{\frac1n:n\in\Bbb N\}\cup\{0\}$. Then $A$ is compact in $\Bbb C$
but the projective cone 
$\Cal H(A\subseteq \Bbb C,\Bbb C) \to  \Cal H(\{z\}\subseteq\Bbb C,\Bbb C)$ 
($z\in A$) does not generate the bornology.
\endproclaim
\demo{Proof}
Let $\Cal B\subseteq \Cal H(A\subseteq \Bbb C,\Bbb C)$ be the set of germs 
of the following locally constant functions $f_n:\{x+iy\in\Bbb C:x\neq 
r_n\} \to  \Bbb C$, with $f_n(x+iy)$ being equal to 0 for $x<r_n$ and 
being equal to 1 for $x>r_n$, where $r_n:=
 \frac2{2n+1}$, for $n\in \Bbb N$.
Then $\Cal B\subseteq \Cal H(A\subseteq \Bbb C,\Bbb C)$ is not bounded, 
otherwise there would exist a neighborhood $W$ of $A$ such that the germ 
of $f_n$ extends to a holomorphic mapping on $W$ for all $n$.
Since every $f_n$ is 0 on some neighborhood of 0, these extensions have 
to be zero on the component of $W$ containing 0, which is not possible, 
since $f_n(\frac1n)=1$.

But on the other hand the set 
$\Cal B_z\subseteq \Cal H(\{z\}\subseteq\Bbb C,\Bbb C)$
of germs at z of all 
germs in $\Cal B$ is bounded, since it contains only the germs of the 
constant functions 0 and 1.
\qed\enddemo

\subheading{\nmb.{3.10}. Spaces of germs of real-analytic functions} 

Let $M$ be a real analytic finite
dimensional manifold. 
If $f:M \to M'$ is a mapping between two such 
manifolds, then $f$ is real analytic if and only if $f$ maps
smooth curves into smooth ones and real analytic curves into
real analytic ones, by \nmb!{2.4}. 

For each real analytic manifold $M$ of real dimension $m$ there is a
\idx{\it complex manifold 
$M_{\Bbb C}$} of complex dimension $m$ containing $M$ as a real analytic 
closed
submanifold, whose germ along $M$ is unique
(\cit!{32}, Proposition 1),  and which can be
chosen even to be a Stein manifold, 
see \cit!{7}, section 3. The complex charts are just
extensions of the real analytic charts of an atlas of $M$.

Real analytic mappings $f:M\to M'$ are the germs
along $M$ of holomorphic mappings $W \to M'_{\Bbb C}$ for open
neighborhoods $W$ of $M$ in $M_{\Bbb C}$.

Let $C^\omega(M,F)$ be the \idx{\it space of real analytic functions}
$f:M\to F$, for any convenient vector space $F$, and let 
$\Cal H(M,\Bbb C) := \Cal H(M\subseteq M_{\Bbb C},\Bbb C)$.
Furthermore, for a subset $A\subseteq M$ let $C^\om(A\subseteq M, \Bbb R)$ 
denotes the \idx{\it space of
germs of real analytic functions} defined near $A$.

\proclaim{\nmb.{3.11}. Lemma} For any subset $A$ of $M$ the 
complexification of the real vector space $C^\om(A\subseteq M,\Bbb R)$
is the complex vector space $\Cal H(A\subseteq M_{\Bbb C},\Bbb C)$.
\endproclaim

\subheading{Definition} For any $A\subseteq M$ of a real analytic
manifold $M$ we will
topologize $C^\om(A\subseteq M,\Bbb R)$ 
as subspace of $\Cal H(A\subseteq M_{\Bbb C},\Bbb C)$, in fact as
the real part of it.

\demo{Proof} Let $f,g\in C^\om(A\subseteq M,\Bbb R)$. They are germs of real 
analytic mappings defined on some open neighborhood of $A$ in $M$.
Inserting complex numbers into the locally convergent Taylor
series in local coordinates
shows, that $f$ and $g$ can be considered as holomorphic mappings 
from some neighborhood $W$ of $A$ in $M_{\Bbb C}$, which have real values
if restricted to $W\cap M$. The mapping $h:=f+ig:W\to \Bbb C$ gives then an 
element of $\Cal H(A\subseteq M_{\Bbb C},\Bbb C)$.

Conversely let $h\in\Cal H(A\subseteq M_{\Bbb C},\Bbb C)$. Then
$h$ is the germ of a holomorphic 
mapping $\widetilde h:W\to \Bbb C$ for some open neighborhood $W$ of
$A$ in $M_{\Bbb C}$.
The decomposition of $h$ into real and imaginary part $f={1\over
2}(h+\bar h)$ and
$g={1\over 2}(h-\bar h)$, which are real analytic maps if restricted to 
$W\cap M$, gives elements of $C^\om(A\subseteq M,\Bbb R)$.

That these correspondences are inverse to each other follows from the 
fact that a holomorphic germ is determined by its restriction to
a germ of mappings $M\supseteq A\to \Bbb C$.\qed
\enddemo

\proclaim{\nmb.{3.12}. Lemma} The inclusion 
$C^\om(M,\Bbb R) \to C^\infty(M,\Bbb R)$
is continuous.
\endproclaim
\demo{Proof} Consider the following diagram, where $W$ is an
open neighborhood of $M$ in $M_{\Bbb C}$.
$$\CD 
C^\om(M,\Bbb R) @>{\text{inclusion}}>> C^\infty(M,\Bbb R)\\
@V\text{direct summand}VV                @VV\text{direct summand}V\\
\Cal H(M\subseteq M_{\Bbb C},\Bbb C) @>\text{inclusion}>> C^\infty(M,\Bbb R^2)\\
@A\text{restriction}AA                   @AA\text{restriction}A\\
\Cal H(W,\Bbb C)      @>\text{inclusion}>>   C^\infty(W,\Bbb R^2)\qed
\endCD$$
\enddemo

\medpagebreak
\proclaim{\nmb.{3.13}. Theorem (Structure of $C^\om(A\subseteq M,\Bbb R)$ 
for closed subsets $A$ of real analytic manifolds $M$)}
The inductive cone
$$C^\om(A\subseteq M,\Bbb R) \gets 
\{\;C^\om(W,\Bbb R): A\subseteq W \underset\text{open}\to\subseteq M\}$$
is regular, i.e. every bounded set is contained and bounded in
some step.

The projective cone
$$C^\om(A\subseteq M,\Bbb R) \to \{\;C^\om(K\subseteq M,\Bbb R):
K\text{ compact in }A \}$$
generates the bornology of $C^\om(A\subseteq M,\Bbb R)$.

If $A$ is even a smooth submanifold, then the following projective cone 
also generates the bornology.
$$C^\om(A\subseteq M,\Bbb R) \to 
\{\;C^\om(\{x\}\subseteq M,\Bbb R):x\in A \}$$

The space $C^\om(\{0\}\subseteq \Bbb R^m,\Bbb R)$ is also the regular 
inductive limit of the spaces 
$\ell^p_r (r\in\Bbb R^m_+)$ for all $1\leq p \leq \infty$.

For general closed $A\subseteq N$ the space $C^\om(A\subseteq M,\Bbb R)$ is 
Montel (hence quasi-complete and reflexive), and ultra-bornological 
(hence a convenient vector space). It is also webbed and conuclear.
If $A$ is compact then it is even a strongly nuclear Silva space and 
its dual is a strongly nuclear Fr\'echet space. 
It is however not a Baire space.
\endproclaim

\demo{Proof} This follows using \nmb!{3.11} from \nmb!{3.4}, 
\nmb!{3.6}, and \nmb!{3.8} by passing to the real parts and from the 
fact that all properties are inherited by complemented subspaces as
$C^\om(A\subseteq M,\Bbb R)$ of $\Cal H(A\subseteq M_{\Bbb C},\Bbb C)$.
\qed\enddemo

\proclaim{\nmb.{3.14} Corollary} A subset $\Cal B\subseteq 
C^\om(\{0\}\subseteq \Bbb R^m,\Bbb R)$ is bounded if and only if
$\Cal B^{(\al)}_0:=\{f^{(\al)}(0):f\in \Cal B\}$ is bounded in
$\Bbb R$ for all  
$\al\in \Bbb N_0^m$ and the poly-radius of convergence for $f\in \Bbb B$ is 
bounded from below by some $r_0\in\Bbb N_0^m$ (or equivalently there 
exists an
$r > 0$ such that 
$\{\tfrac {f^{(\al)}}{\al!}r^{|\al|}:f\in\Bbb B, \al\in \Bbb N_0^m\}$ is 
bounded in $\Bbb R$).
\endproclaim
\demo{Proof} The space $C^\om(\{0\}\subseteq \Bbb R^m,\Bbb R)$ is the 
regular inductive limit of the spaces 
$\ell^p (r\in\Bbb R^m_+)$ for $p$ equal to 1 or $\infty$ by \nmb!{3.13}. 
Hence $\Cal B$ is bounded if and only if it is contained and bounded in 
$\ell^p_r$ for some
$r\in\Bbb R^m_+$. This shows the equivalence with the first 
condition using $p=1$ and the equivalence with the second condition using 
$p=\infty$.
\qed\enddemo

\heading\totoc \nmb0{4}. A uniform boundedness principle \endheading

\proclaim{\nmb.{4.1}. Lemma} Let $(E,E')$ be a dual pairing and let
$\Cal S$ be a point separating set of bounded linear mappings with 
common domain $(E,E')$. Then the following conditions are equivalent.
\roster
\item If $F$ is a Banach space (or even a complete dual pairing $(F,F')$)
and $f:F \to  E$ is linear and $\la \o f$ 
is bounded for all $\la \in \Cal S$, then f is bounded.
\item If $B\subseteq E$ is absolutely convex such that  $\la(B)$
is bounded for all $\la \in \Cal S$ and the normed space $E_B$
generated by $B$ is complete, then $B$ is bounded in $E$.
\item Let $(b_n)$ be an unbounded sequence in $E$ with
$\la(b_n)$ bounded for all $\la \in \Cal S$, then there is
some $(t_n) \in \ell^1$ such that $\sum t_n\,b_n$ does not
converge in $E$ for the weak topology induced by $\Cal S$.
\endroster
\endproclaim

{\bf Definition.} We say that $(E,E')$ satisfies the {\it uniform
$\Cal S$-boundedness principle} \idx{\ign{uniform S-boundedness principle}} 
if these equivalent conditions are satisfied. 

\demo{Proof}
\therosteritem 1 $\Rightarrow$ \therosteritem 3 :
Suppose that \therosteritem 3 is not satisfied. 
So let $(b_n)$ be an unbounded sequence in $E$ such that
$\la(b_n)$ is bounded for all $\la \in \Cal S$, and such that for all
$(t_n) \in \ell^1$ the series $\sum t_n\,b_n$
converges in $E$ for the weak topology induced by $\Cal S$.
We define a linear mapping $f:\ell^1 \to  E$ by $f(t_n) = \sum t_n\,b_n$.
It is easy to check that $\la \o f$ is bounded, hence by
\therosteritem 1 the image of the closed unit ball, which
contains all $b_n$, is bounded. Contradiction.

\therosteritem{3} $\Rightarrow$ \therosteritem{2}:
Let $B \subseteq E$ be absolutely convex such that  $\la(B)$
is bounded for all $\la \in \Cal S$ and that the normed space $E_B$
generated by $B$ is complete, and suppose that $B$ is unbounded.
Then $B$ contains an unbounded sequence $(b_n)$, so by
\therosteritem{3} there is some $(t_n) \in \ell^1$ such that
$\sum t_n\,b_n$ does not converge in $E$ for the weak topology
induced by $\Cal S$. But $\sum t_n\,b_n$ is easily seen to
be a Cauchy sequence in $E_B$ and thus converges even bornologically, a 
contradiction.

\therosteritem{2} $\Rightarrow$ \therosteritem{1}:
Let the bornology of $F$ be complete, and let $f:F \to  E$ be linear such that 
$\la \o f$ is bounded for all $\la \in \Cal S$. It suffices to show
that $f(B)$, the image of an absolutely convex bounded set $B$
in $F$ with $F_B$ complete,
is bounded. Then $\la(f(B))$ is bounded for all $\la \in \Cal S$, 
the normed space
$E_{f(B)}$ is a quotient of $F_B$, hence complete. By
\therosteritem{2} the set $f(B)$ is bounded.
\qed\enddemo

\proclaim{\nmb.{4.2}. Lemma} A complete dual pair $(E,E')$
satisfies the uniform $\Cal S$-boundedness principle for each
point separating set $\Cal S$ of bounded linear mappings on $E$ 
if and only if there exists no strictly weaker ultrabornological topology than
the natural bornological topology of $(E,E')$.
\endproclaim
\demo{Proof} $(\Rightarrow)$ Let $\ta$ be an ultrabornological
topology on $E$ which is strictly weaker than the natural
bornological topology. Since every ultra-bornological space is an 
inductive limit of Banach spaces, cf. \cit!{10}, 13.1.2, 
there exists a Banach space
$F$ and a continuous linear mapping $f:F\to (E,\ta)$ which is not
continuous into $E$. Let $\Cal S=\{Id:E\to(E,\tau)\}$.
Now $f$ does not satisfy \nmb!{4.1}.\therosteritem1.

$(\Leftarrow)$ If $\Cal S$ is a point separating set of bounded
linear mappings, the ultrabornological topology given by the
inductive limit of the spaces $E_B$ with $B$ satisfying
\nmb!{4.1}.\therosteritem2 equals the natural bornological
topology of $(E,E')$.
Hence \nmb!{4.1}.\therosteritem2 is satisfied.
\qed\enddemo

\proclaim{\nmb.{4.3}. Lemma} Let $\Cal F$ be a set of bounded linear 
mappings $f:E\to E_f$ between dual pairings, let $\Cal S_f$ 
be a point separating set of bounded linear mappings on $E_f$ for every 
$f\in \Cal F$, and let $\Cal S:=\bigcup_{f\in \Cal F}f^*(\Cal S_f) =
\{g\o f: f\in \Cal F, g\in \Cal S_f\}$. If $\Cal F$ generates the 
bornology and $E_f$ satisfies
the uniform $\Cal S_f$-boundedness principle for all $f\in \Cal F$, then 
$E$ satisfies the uniform $\Cal S$-boundedness principle.
\endproclaim
\demo{Proof}
We check the condition (1) of \nmb!{4.1}. So assume $h:F\to E$ is a linear 
mapping for which $g\o f\o h$ is bounded for all $f\in \Cal F$ and 
$g\in \Cal S_f$. Then $f\o h$ is bounded by the uniform $\Cal S_f$-
boundedness principle for $E_f$. Consequently $h$ is bounded since $\Cal F$ 
generates the bornology of $E$.
\qed\enddemo

\proclaim{\nmb.{4.4}. Theorem} A locally convex space which is webbed
satisfies the uniform $\Cal S$-boun\-ded\-ness principle for any
point separating set of bounded linear functionals.
\endproclaim
\demo{Proof} Since the bornologification of a webbed space is webbed, 
cf. \cit!{10}, 13.3.3 and 13.3.1, we may assume that $E$ is 
bornological, and hence that every bounded linear functional is 
continuous, cf\. \cit!{10}, 13.3.1.
Now the closed graph principle, cf. \cit!{10}, 56.4.1, applies 
to any mapping satisfying the assumptions of \nmb!{4.1}.1. 
\qed\enddemo

\proclaim{\nmb.{4.5}. Theorem (Holomorphic uniform boundedness 
principle)} \newline
For any closed subset $A \subseteq N$ of a complex manifold $N$ the
locally convex space $\Cal H(A\subseteq N,\Bbb C)$ satisfies the uniform 
$\Cal S$-boundedness principle for every point separating set
$\Cal S$ of bounded linear functionals.
\endproclaim
\demo{Proof} This is a immediate consequence of \nmb!{4.4} and 
\nmb!{3.6}.
\qed\enddemo

\demo{Direct proof of a particular case} 
We prove the theorem for a closed smooth submanifold $A\subseteq \Bbb C$ and
the set $\Cal S$ of all iterated derivatives at points in $A$.

Let us suppose first
that $A$ is the point 0.
We will show that condition \nmb!{4.1}.3 is satisfied. Let
$(b_n)$ be an unbounded sequence in $\Cal H(\{0\},\Bbb C)$ such
that each Taylor coefficient $b_{n,k} =
\frac1{k!}\,b_n^{(k)}(0)$  is bounded with respect to $n$:
$$\sup\{\, |b_{n,k}|:n\in \Bbb N\,\} < \infty. \tag1$$
We have to find $(t_n) \in \ell^1$ such that $\sum_n t_n\,b_n$ is
no longer the germ of a holomorphic function at $0$.

Each $b_n$ has positive radius of convergence, in particular
there is an $r_n >0$ such that 
$$\sup\{\,|b_{n,k}\,r_n^k|:k\in\Bbb N\,\} < \infty. \tag2$$
By theorem \nmb!{3.4} the space $\Cal H(\{0\},\Bbb C)$ is a regular
inductive limit of spaces $\ell^\infty_r$. Hence a subset $\Cal B$
is bounded in $\Cal H(\{0\},\Bbb C)$ if and only if there
exists an $r>0$ such that  
$$\left\{\,\frac 1{k!}b^{(k)}(0)\,r^k: b\in \Cal B, k 
\in \Bbb N\,\right\}$$
is bounded.
That the sequence $(b_n)$ is unbounded thus means that for all
$r>0$ there are $n$ and $k$ such that
$|b_{n,k}| > (\tsize\frac 1r)^k$.
We can even choose $k>0$ for otherwise the set
$\{\,b_{n,k}r^k:n,k\in \Bbb N, k>0\,\}$ is bounded, so only 
$\{\,b_{n,0}:n\in \Bbb N\,\}$ can be unbounded. This
contradicts \thetag1.

Hence for each $m$ there are $k_m>0$ such that 
$\Cal N_m := \{\,n \in \Bbb N: |b_{n,k_m}| > m^{k_m}\,\}$  
is not empty.
We can choose $(k_m)$ strictly increasing, for if they were
bounded, $|b_{n,k_m}| < C$ for some $C$ and all $n$ by \thetag1, but 
$|b_{n_m,k_m}| > m^{k_m} \to  \infty$ for some $n_m$.

Since by \thetag1 the set $\{\,b_{n,k_m}: n\in \Bbb N\,\}$ is
bounded, we can choose $n_m \in \Cal N_m$ such that 
$$\aligned  
|b_{n_m,k_m}|&\geq {\tsize\frac12} |b_{j,k_m}| \quad \text{ for } j>n_m\\
|b_{n_m,k_m}| &> m^{k_m}\endaligned \tag3$$

We can choose also $(n_m)$ strictly increasing, for if they were
bounded we would get $|b_{n_m,k_m}r^{k_m}| < C$ for
some $r>0$ and $C$ by \thetag 2. But $(\frac 1m)^{k_m} \to  0$. 

We pass now to the subsequence $(b_{n_m})$ which we denote again
by $(b_m)$. We put
$$t_m := \operatorname{sign}\left(\frac 1{b_{m,k_m}}
\sum_{j<m}t_j\,b_{j,k_m}\right)\cdot \frac 1{4^m}.\tag4$$
Assume now that $b_\infty = \sum_m t_m\,b_m$ converges weakly
with respect to $\Cal S$ to a holomorphic germ. Then
its Taylor series is 
$b_\infty(z) = \sum_{k\geq 0} b_{\infty,k}\,z^k$,
where the coefficients are given by 
$b_{\infty,k}=\sum_{m\geq 0} t_m\,b_{m,k}$.
But we may compute as follows, using \thetag3 and \thetag4 :
$$\allowdisplaybreaks\align
|b_{\infty,k_m}| &\geq 
\biggl|\sum_{j\leq m}t_j\,b_{j,k_m}\biggr| - 
\sum_{j>m}|t_j\,b_{j,k_m}| =\\
&\aligned = \biggl|\sum_{j<m}t_j\,b_{j,k_m}\biggr| &+ 
|t_m\,b_{m,k_m}| \qquad\text{(same sign)}\\
& -\sum_{j>m}|t_j\,b_{j,k_m}| \geq\endaligned\\
&\geq 0 + |b_{m,k_m}|\cdot\left(|t_m| - 2\sum_{j>m}|t_j|\right)=\\
&= |b_{m,k_m}|\cdot \frac 1{3\cdot 4^m} \geq
\frac{m^{k_m}}{3\cdot 4^m}.
\endalign$$
So $|b_{\infty,k_m}|^{1/{k_m}}$ goes to $\infty$, hence
$b_\infty$ cannot have a positive radius of convergence, a
contradiction. So the theorem follows for the space $\Cal
H(\{t\},\Bbb C)$.

Let us consider now an arbitrary closed smooth submanifold 
$A \subseteq \Bbb C$. 
By \nmb!{3.8} the projective cone 
$\Cal H(A\subseteq N,\Bbb C) \to  \{\Cal H(\{z\}\subseteq N,\Bbb C), z\in A\}$ 
generates the bornology. Hence the result follows from 
the case where $A=\{0\}$ by \nmb!{4.3}.\qed\enddemo

\proclaim{\nmb.{4.6}. Theorem (Special real analytic uniform boundedness
principle)} 
For any closed subset $A\subseteq M$ of a real
analytic manifold $M$, the space $C^\om(A\subseteq M,\Bbb R)$
satisfies the uniform $\Cal S$-boundedness principle for any
point separating set $\Cal S$ of bounded linear functionals. 

If $A$ has no isolated points and $M$ is 1-dimensional this applies 
to the set of all
point evaluations $\ev_t$, $t\in A$. 
\endproclaim

\demo{Proof}
Again this follows from \nmb!{4.4} using now
\nmb!{3.13}. If $A$ has no isolated points and $M$ is 1-dimensional
the point evaluations are separating, by the uniqueness theorem for 
holomorphic functions.
\qed\enddemo

\demo{Direct proof of a particular case} We show that
$C^\om(\Bbb R,\Bbb R)$ satisfies the uniform $\Cal
S$-boun\-ded\-ness principle for the set $\Cal S$ of all point evaluations. 

We check property \nmb!{4.1}.2. Let 
$\Cal B \subseteq C^\om(\Bbb R,\Bbb R)$ be absolutely convex such
that $\ev_t(\Cal B)$ is bounded for all $t$ and such that 
$C^\om(\Bbb R,\Bbb R)_B$ is complete. We have to show that $\Cal B$
is complete.

By lemma \nmb!{3.12} the set $\Cal B$ satisfies the conditions of 
\nmb!{4.1}.2 in
the space $C^\infty(\Bbb R,\Bbb R)$. Since 
$C^\infty(\Bbb R,\Bbb R)$ satisfies the uniform $\Cal
S$-boundedness principle, cf. \cit!{5}, the set 
$\Cal B$ is
bounded in $C^\infty(\Bbb R,\Bbb R)$. Hence all iterated
derivatives at points are bounded on $\Cal B$, and a fortiori the
conditions of \nmb!{4.1}.2 are satisfied for $\Cal B$ in $\Cal H(\Bbb
R,\Bbb C)$. By the particular case of theorem \nmb!{4.5} the set
$\Cal B$ is bounded in  
$\Cal H(\Bbb R,\Bbb C)$ and hence also in the direct summand 
$C^\om(\Bbb R,\Bbb R)$.\qed\enddemo

\heading\totoc \nmb0{5}. Cartesian closedness \endheading

\proclaim{\nmb.{5.1}. Theorem} The real analytic curves 
in $C^\om(\Bbb R,\Bbb R)$ correspond exactly to the real
analytic functions $\Bbb R^2 \to  \Bbb R$.
\endproclaim
\demo{Proof} ($\Rightarrow$) Let $f:\Bbb R \to  C^\om(\Bbb
R,\Bbb R)$ be a real analytic curve. Then
$f: \Bbb R \to  C^\om(\{t\},\Bbb R)$ is also real analytic. 
We use theorems \nmb!{3.13} and \nmb!{1.6} to conclude that $f$ is
even a topologically real analytic curve in $C^\om(\{t\},\Bbb R)$.
By lemma \nmb!{1.7} for every $s\in \Bbb R$ the curve $f$ can
be extended to a holomorphic mapping from an open
neighborhood of $s$ in $\Bbb C$ to the
complexification (\nmb!{3.11}) $\Cal H(\{t\},\Bbb C)$ of 
$C^\om(\{t\},\Bbb R)$.

 From \nmb!{3.4} it follows that $\Cal H(\{t\},\Bbb C)$ is the regular
inductive limit of all spaces $\Cal H(U,\Bbb C)$, where $U$
runs through some neighborhood basis of $t$ in $\Bbb C$.
Lemma \nmb!{1.8} shows that $f$ is a holomorphic mapping $V\to  \Cal
H(U,\Bbb C)$ for some open neighborhoods $U$ of $t$
and $V$ of $s$ in $\Bbb C$.

By the exponential law for
holomorphic mappings (see \nmb!{1.3}) the canonically associated
mapping $f\sphat: V\times U \to  \Bbb C$ is holomorphic.
So its restriction is a real analytic function $\Bbb R\times \Bbb
R \to \Bbb R$ near $(s,t)$.

($\Leftarrow$) Let $f:\Bbb R^2 \to  \Bbb R$ be a real analytic mapping.
Then $f(t,\quad)$ is real analytic, so the associated mapping 
$f\spcheck:\Bbb R \to  C^\om(\Bbb R,\Bbb R)$ makes sense. It
remains to show that it is real analytic. Since the mappings $C^\om(\Bbb R,\Bbb R) 
\to  C^\om(K,\Bbb R)$ generate the bornology, by \nmb!{3.13}, it is by 
\nmb!{1.11} enough 
to show that $f\spcheck:\Bbb R \to  C^\om(K,\Bbb R)$ is real analytic 
for each compact $K\subseteq \Bbb R$, which may be checked
locally near each $s\in \Bbb R$.

$f:\Bbb R^2 \to  \Bbb R$ extends to a holomorphic function on an
open neighborhood $V\times U$ of $\{s\}\times K$ in $\Bbb C^2$.
By cartesian closedness for the holomorphic setting
the associated mapping 
$f\spcheck:V\to  \Cal H(U,\Bbb C)$ is holomorphic, so its
restriction $V\cap\Bbb R \to  C^\om(U\cap \Bbb R,\Bbb R)\to 
C^\om(K,\Bbb R)$ is real analytic as required.
\qed\enddemo

\subheading{\nmb.{5.2}. Remark} From \nmb!{5.1} it follows that the curve
$c:\Bbb R\to C^\om(\Bbb R,\Bbb R)$ defined in \nmb!{1.1} is real analytic, 
but it is not topologically real analytic.
In particular, it does not factor locally
to a real analytic curve into  some Banach space 
$C^\om(\Bbb R,\Bbb R)_B$ for a bounded subset $B$ and it has no 
holomorphic extension to a mapping defined on a neighborhood of $\Bbb R$
in $\Bbb C$ with values in the complexification $\Cal H(\Bbb R,\Bbb C)$
of $C^\om(\Bbb R,\Bbb R)$, cf. \nmb!{1.7}.

\proclaim{\nmb.{5.3}. Lemma} For a real analytic manifold $M$, 
the bornology on the space $C^\omega(M,\Bbb R)$ is induced by the
following cone. 
$$C^\omega(M,\Bbb R) @>{c^*}>> C^\al(\Bbb R,\Bbb R)$$
for all $C^\al$-curves
$c:\Bbb R\to M$, where $\al$ equals $\infty$ and $\om$.
\endproclaim
\demo{Proof} The maps $c^*$ are bornological since $C^\om(M,\Bbb
R)$ is convenient by \nmb!{3.13}, and by the uniform $\Cal
S$-boundedness principle \nmb!{4.6} for $C^\om(\Bbb R,\Bbb R)$ and by 
\cit!{5}, 4.4.7 for $C^\infty(\Bbb R,\Bbb R)$ 
it suffices to check that 
$\ev_t\o c^* = \ev_{c(t)}$ is bornological, which is obvious.

Conversely we consider the identity mapping $i$ from the space
$E$ into $C^\om(M,\Bbb R)$, where $E$ is the vector space
$C^\om(M,\Bbb R)$, but with the locally convex structure induced
by the cone. \newline
{\bf Claim.} The bornology of $E$ is complete.\newline
The spaces $C^\om(\Bbb R,\Bbb R)$ and $C^\infty(\Bbb R,\Bbb R)$
are convenient by \nmb!{3.13} and \nmb!{1.3}, respectively. 
So their product 
$$\prod_cC^\om(\Bbb R,\Bbb R)\x\prod_cC^\infty(\Bbb R,\Bbb R)$$ 
is also convenient. By
theorem \nmb!{2.4},\therosteritem1 $\Leftrightarrow$
\therosteritem5, the embedding of $E$ into this product has 
closed image, hence the bornology of $E$ is complete.

Now we may apply the uniform $S$-boundedness principle
for $C^\om(M,\Bbb R)$ (\nmb!{4.6}), since obviously 
$\ev_p\o i = \ev_0\o c_p^*$ is bounded, where $c_p$ is the
constant curve with value $p$, for all $p\in M$. \qed\enddemo

\subheading{\nmb.{5.4}. Structure on $C^\om(U,F)$} Let $(E,E')$ be
a dual pair of real vector spaces and let $U$ be $c^\infty$-open in E.
We equip the space
$C^\om(U,\Bbb R)$ of all real analytic functions (cf.
\nmb!{2.6}) with the dual space consisting of all linear functionals
induced from the families of mappings
$$\gather C^\om(U,\Bbb R) @>{c^*}>> C^\om(\Bbb R,\Bbb R),\text{
for all }c\in C^\om(\Bbb R,U)\\ 
C^\om(U,\Bbb R) @>{c^*}>> C^\infty(\Bbb R,\Bbb R),\text{
for all }c\in C^\infty(\Bbb R,U).
\endgather$$
For a finite dimensional vector spaces $E$ this definition gives
the same bornology as the one defined in \nmb!{3.10}, by lemma 
\nmb!{5.3}.

If $(F,F')$ is another dual pair, we equip the space
$C^\om(U,F)$ of all real analytic mappings (cf. \nmb!{2.6}) with
the dual induced by the family of mappings 
$$ C^\om(U,F) @>{\la_*}>> C^\om(U,\Bbb R),\text{ for all }\la\in
F'.$$

\proclaim{\nmb.{5.5}. Lemma} Let $(E,E')$ and $(F,F')$ be complete dual
pairs and let $U\subseteq E$ be $c^\infty$-open. Then $C^\om(U,F)$ is
complete.\endproclaim

\demo{Proof}
This follows immediately from the fact that $C^\om(U,F)$ can be considered
as closed subspace of the product of factors $C^\om(U,\Bbb R)$
indexed by all $\la\in F'$. And $C^\om(U,\Bbb R)$ can be considered as
closed subspace of the product of the factors $C^\om(\Bbb R,\Bbb R)$ indexed 
by all $c\in C^\om(\Bbb R,U)$ and the factors $C^\infty(\Bbb R,\Bbb R)$
indexed by all $c\in C^\infty(\Bbb R,U)$. Since all factors are complete
so are the closed subspaces.\qed
\enddemo

\proclaim{\nmb.{5.6}. Lemma (General real analytic uniform
boundedness principle)} Let $E$ and $F$ be convenient 
vector spaces and $U\subseteq E$ be $c^\infty$-open. Then 
$C^\om(U,F)$ satisfies the uniform $\Cal S$-boundedness principle, 
where $\Cal S:=\{ev_x:x\in U\}$.
\endproclaim
\demo{Proof} The complete bornology of $C^\om(U,F)$ is by 
definition induced by the maps $c^*:C^\om(U,F)\to C^\om(\Bbb R,F)$ 
($c\in C^\om(\Bbb R,U)$) together with the maps 
$c^*:C^\om(U,F)\to C^\infty(\Bbb R,F)$ 
($c\in C^\infty(\Bbb R,U)$). Both spaces $C^\om(\Bbb R,F)$ and 
$C^\infty(\Bbb R,F)$ satisfy the uniform $\Cal T$-boun\-ded\-ness 
principle, where $\Cal T:=\{ev_t:t\in \Bbb R\}$, by \nmb!{4.6} and
\cit!{5}, 4.4.7, respectively.
Hence $C^\om(U,F)$ satisfies the uniform $\Cal S$-boundedness
principle by lemma \nmb!{4.3}, 
since $ev_t\o c^*=ev_{c(t)}$.
\qed\enddemo

\subheading{\nmb.{5.7}. Definition} Let $(E,E')$ and $(F,F')$ be 
complete dual pairs.
We denote by $L(E,F)$ the \idx{\it space of linear real analytic
mappings} from E to F,  
which are by \nmb!{1.9} exactly the bounded linear mappings. 
Furthermore, if $E$ and $F$ are convenient vector spaces, these are 
exactly the morphisms in the sense of dual pairs, since $f$ is
bounded if and only if $\la\o 
f\in E^b$ for all $\la\in F^b$.

\proclaim{\nmb.{5.8}. Lemma (Structure on $L(E,F)$)}
The following structures on $L(E,F)$ are the same:
\roster
\item The bornology of pointwise boundedness, i.e. the bornology 
      induced by the cone $(ev_x:L(E,F)\to F, x\in E)$.
\item The bornology of uniform boundedness on bounded sets in 
      $E$, i.e. a set $\Cal B\subseteq L(E,F)$ is bounded if and only if $\Cal 
      B(B)\subseteq F$ is bounded for every bounded $B\subseteq E$.
\item The bornology induced by the inclusion $L(E,F) \to 
      C^\infty(E,F)$.
\item The bornology induced by the inclusion $L(E,F) \to 
      C^\om(E,F)$.
\endroster
\endproclaim

The space $L(E,F)$ will from now on be the convenient vector 
space having as structure that described in the previous lemma.
Thus $L(E,F)$ is a convenient vector space, 
by \cit!{5}, 3.6.3. 
In particular this is true for $E'=L(E,\Bbb R)$.

So a mapping $f$ into $L(E,F)$ is real analytic if and only if
the composites $ev_x\o f$ are real analytic for all $x\in E$, by
\nmb!{1.11}.

\demo{Proof} That the bornology in \therosteritem 1, 
\therosteritem 2 and \therosteritem 3 are the same was shown in 
\cit!{5}, 3.6.4 and 4.4.24. Since
$C^\om(E,F)\to C^\infty(E,F)$  
is continuous by definition of the structure on $C^\om(E,F)$ the 
bornology in \therosteritem 4 is finer than that in 
\therosteritem 1. The bornology given in \therosteritem 4 is 
complete, since the point-evaluations $ev_x:C^\om(E,F)\to F$ are 
continuous, and linearity of a mapping $E\to F$ can be checked by applying 
them. Furthermore
$L(E,F)$ with the bornology given in \therosteritem 4 satisfies 
the uniform $\Cal S$-boundedness theorem, since $C^\om(E,F)$ 
does, by \nmb!{5.6}. So the 
identity on $L(E,F)$ with the bornology given in \therosteritem 1 
to that given in \therosteritem 4 is bounded. 
\qed\enddemo

The following two results 
will be generalized in \nmb!{6.3}.
At the moment we will make use of the following lemma only in case where 
$E=C^\infty(\Bbb R,\Bbb R)$.

\proclaim{\nmb.{5.9}. Lemma} $L(E,C^\om(\Bbb R,\Bbb R)) \cong 
C^\om(\Bbb R,E')$ as vector spaces, for any convenient
vector space $E$.
\endproclaim
\demo{Proof} For $c\in C^\om(\Bbb R,E')$ consider $\tilde c(x)
:= \ev_x\o c\in C^\om(\Bbb R,\Bbb R)$ for $x\in E$. By the
uniform $\Cal S$-boundedness principle \nmb!{4.6} for $\Cal S=\{\ev_t:t\in
\Bbb R\}$ the linear mapping $\tilde c$ is bounded, since 
$\ev_t\o\tilde c=c(t)\in E'$.

If conversely $\ell\in L(E,C^\om(\Bbb R,\Bbb R))$, we consider
$\tilde \ell(t) = \ev_t\o \ell \in E':=L(E,\Bbb R)$ for $t\in \Bbb R$. 
Since the bornology of $E'$ is generated by $\Cal S:=\{ev_x:x\in E\}$, 
$\tilde \ell:\Bbb R \to  E'$ is real analytic, for 
$\ev_x\o\tilde \ell=\ell(x)\in C^\om(\Bbb R,\Bbb R)$.
\qed\enddemo

\proclaim{\nmb.{5.10}. Corollary} We have 
$C^\infty(\Bbb R,C^\om(\Bbb R,\Bbb R)) \cong 
C^\om(\Bbb R,C^\infty(\Bbb R,\Bbb R))$ as vector spaces.
\endproclaim
\demo{Proof} $ C^\infty(\Bbb R,\Bbb R)'$ is the free convenient
vector space over $\Bbb R$ by \cit!{5}, 5.1.8, 
and $C^\om(\Bbb R,\Bbb R)$ is convenient, we have
$$\align C^\infty(\Bbb R,C^\om(\Bbb R,\Bbb R)) &\cong
    L(C^\infty(\Bbb R,\Bbb R)',C^\om(\Bbb R,\Bbb R))\\
&\cong C^\om(\Bbb R,C^\infty(\Bbb R,\Bbb R)'')\qquad
\text{by lemma \nmb!{5.9}}\\
&\cong C^\om(\Bbb R,C^\infty(\Bbb R,\Bbb R)),
\endalign$$
by reflexivity of $C^\infty(\Bbb R,\Bbb R)$, see 
\cit!{5}, 5.4.16.
\qed\enddemo

\proclaim{\nmb.{5.11}. Theorem}
Let $(E,E')$ be a complete dual pair, let $U$ be $c^\infty$-open
in $E$, let $f:\Bbb R\x U\to  \Bbb R$ be
a real analytic mapping and let $c\in C^\infty(\Bbb R,U)$. Then
$c^*\o\check f:\Bbb R \to  C^\om(U,\Bbb R) \to  C^\infty(\Bbb R,\Bbb
R)$ is real analytic.
\endproclaim

This result on the mixing of $C^\infty$ and $C^\om$ will become quite
essential in the proof of cartesian closedness. It will be generalized 
in \nmb!{6.4}, see also \nmb!{8.9} and \nmb!{8.14}.

\demo{Proof} Let $I\subseteq\Bbb R$ be open and relatively compact, let
$t\in\Bbb R$ and $k\in\Bbb N$. Now choose an open and relatively compact
$J\subseteq\Bbb R$ containing the closure $\bar I$ of $I$. There
is a bounded subset $B\subseteq E$ such that
$c\mid J: J \to  E_B$ is a $\Cal Lip^k$-curve
in the Banach space $E_B$ generated
by $B$. This is \cit!{13}, Folgerung on p.114.
Let $U_B$ denote the open subset $U\cap E_B$ of the Banach space $E_B$.
Since the inclusion $E_B \to  E$ is continuous, $f$ is real analytic
as a function $\Bbb R\x U_B \to  \Bbb R \x U \to  \Bbb R$.
Thus by \nmb!{2.4} there is a holomorphic extension
$f:V\x W\to \Bbb C$ of $f$ to an open set
$V\x W\subseteq\Bbb C\x (E_B)_{\Bbb C}$ containing the compact set
$\{t\}\x c(\bar I)$. By cartesian
closedness of the category of holomorphic mappings
$\check f: V\to  \Cal
H(W,\Bbb C)$ is holomorphic.
Now recall that the bornological structure of $\Cal H(W,\Bbb C)$
is induced by that of $C^\infty(W,\Bbb C):=C^\infty(W,\Bbb R^2)$. And
$c^*:C^\infty(W,\Bbb C)\to \Cal Lip^k(I,\Bbb C)$ is a bounded
$\Bbb C$-linear map, by \cit!{5}.
Thus $c^*\o \check f:V\to \Cal Lip^k(I,\Bbb C)$ is holomorphic, and hence
its restriction to $\Bbb R\cap V$, which has values in $\Cal Lip^k(I,\Bbb R)$,
is (even topologically) real analytic by \nmb!{1.7}. Since $t\in
\Bbb R$ was arbitrary we 
conclude that $c^*\o \check f:\Bbb R\to \Cal Lip^k(I,\Bbb R)$ is real analytic. But the
bornology of $C^\infty(\Bbb R,\Bbb R)$ is generated by the inclusions
into $\Cal Lip^k(I,\Bbb R)$, \cit!{5}, 4.2.7, and hence
$c^*\o \check f:\Bbb R \to C^\infty(\Bbb R,\Bbb R)$ is real analytic.
\qed\enddemo

\proclaim{\nmb.{5.12}. Theorem (Cartesian closedness)} The
category of real analytic  
mappings between complete dual pairs of real vector spaces 
is cartesian closed.
More precisely, for complete dual pairs $(E,E')$, $(F,F')$ and $(G,G')$
and $c^\infty$-open sets $U\subseteq E$ and $W\subseteq G$ 
a mapping $f:W\x U\to  F$ is real analytic if and only
if $\check f:W\to  C^\om(U,F)$ is real analytic.
\endproclaim
\demo{Proof}{\bf Step 1.} The theorem is true for $G=F=\Bbb R$.

$(\Leftarrow)$ Let $\check f:\Bbb R\to  C^\om(U,\Bbb
R)$ be $C^\om$. We have to show that $f:\Bbb R\x U\to  \Bbb R$ is $C^\om$.
We consider a curve $c_1:\Bbb R\to \Bbb R$ and a curve $c_2:\Bbb R\to U$.

If the $c_i$ are $C^\infty$, then $c_2^*\o\check f:\Bbb R\to 
C^\om(U,\Bbb R) \to  C^\infty(\Bbb R,\Bbb R)$ is $C^\om$ by
assumption, hence is $C^\infty$, so $c_2^*\o\check f\o c_1:\Bbb
R\to C^\infty(\Bbb R,\Bbb R)$ is $C^\infty$. By cartesian
closedness of smooth mappings, 
$(c_2^*\o\check f\o c_1)^{\wedge }= f\o (c_1\x c_2):\Bbb
R^2\to \Bbb R$ is $C^\infty$. By composing with the diagonal mapping 
$\De:\Bbb R\to \Bbb R^2$ we obtain that $f\o (c_1,c_2):\Bbb R\to \Bbb R$ is
$C^\infty$. 

If the $c_i$ are $C^\om$, then $c_2^*\o\check f:\Bbb R\to 
C^\om(U,\Bbb R) \to  C^\om(\Bbb R,\Bbb R)$ is $C^\om$ by
assumption, so $c_2^*\o\check f\o c_1:\Bbb
R\to C^\om(\Bbb R,\Bbb R)$ is $C^\om$. By theorem \nmb!{5.1}
the associated map
$(c_2^*\o\check f\o c_1)^{\wedge }= f\o (c_1\x c_2):\Bbb
R^2\to \Bbb R$ is $C^\om$. So $f\o (c_1,c_2):\Bbb R\to \Bbb R$ is
$C^\om$. 

$(\Rightarrow)$ Let $f:\Bbb R\x U \to  \Bbb R$ be $C^\om$. We
have to show that $\check f:\Bbb R\to  C^\om(U,\Bbb R)$ is real analytic. 
Obviously $\check f$ has values in this space. 
We consider a curve $c:\Bbb R\to  U$.

If $c$ is $C^\infty$, then by theorem \nmb!{5.11} the
associated mapping $(f\o(Id\x c))\spcheck=c^*\o \check f:
\Bbb R\to  C^\infty(\Bbb R,\Bbb R)$ is $C^\om$.

If $c$ is $C^\om$, then $f\o (Id\x c):\Bbb R\x\Bbb R\to 
\Bbb R\x U \to  \Bbb R$ is $C^\om$. By theorem \nmb!{5.1} the
associated mapping $(f\o(Id\x c))\spcheck=c^*\o \check f:\Bbb
R\to  C^\om(\Bbb R,\Bbb R)$ is $C^\om$.

\smallskip
{\bf Step 2.} The theorem is true for $F=\Bbb R$.

$(\Leftarrow)$ Let $\check f:W\to  C^\om(U,\Bbb R)$ 
be $C^\om$. We have to show that $f:W\x U\to  \Bbb R$ is $C^\om$.
We consider a curve $c_1:\Bbb R\to W$ and a curve $c_2:\Bbb R\to U$.

If the $c_i$ are $C^\infty$, then $c_2^*\o\check f:W \to 
C^\om(U,\Bbb R) \to  C^\infty(\Bbb R,\Bbb R)$ is $C^\om$ by
assumption, hence is $C^\infty$, so $c_2^*\o\check f\o c_1:\Bbb
R\to C^\infty(\Bbb R,\Bbb R)$ is $C^\infty$. By cartesian
closedness of smooth mappings, the associated mapping 
$(c_2^*\o\check f\o c_1)^{\wedge }= f\o (c_1\x c_2):\Bbb
R^2\to \Bbb R$ is $C^\infty$. So $f\o (c_1,c_2):\Bbb R\to \Bbb R$ is
$C^\infty$. 

If the $c_i$ are $C^\om$, then $\check f\o c_1:\Bbb R\to W\to 
C^\om(U,\Bbb R)$ is $C^\om$ by
assumption, so by step 1 the mapping $(\check f\o c_1)^{\wedge}
= f\o (c_1\x Id_U) :\Bbb R\x U\to  \Bbb R$ is $C^\om$. 
Hence 
$$f\o(c_1,c_2)= f\o(c_1\x Id_U)\o(Id,c_2):\Bbb R\to  \Bbb R$$
is $C^\om$.

$(\Rightarrow)$ Let $f:W \x U \to  \Bbb R$ be $C^\om$. We
have to show that $\check f:W \to  C^\om(U,\Bbb R)$ is real analytic.
 Obviously $\check f$ has values in this space. 
We consider a curve $c_1:\Bbb R\to  W$.

If $c_1$ is $C^\infty$, we consider a second curve 
$c_2:\Bbb R\to  U$. If $c_2$ is $C^\infty$, 
then $f\o (c_1\x c_2):\Bbb R\x \Bbb R \to 
W \x U \to  \Bbb R$ is $C^\infty$. By cartesian closedness the
associated mapping 
$$(f\o(c_1 \x c_2))\spcheck = c_2^*\o \check f \o c_1: 
\Bbb R\to  C^\infty(\Bbb R,\Bbb R)$$
is $C^\infty$. If
$c_2$ is $C^\om$, the mapping $f \o (Id_W\x c_2): W\x \Bbb R
\to \Bbb R$ and also the flipped one 
$(f\o(Id_W\x c_2))\sptilde:\Bbb R\x W \to  \Bbb R$
are $C^\om$, hence by theorem \nmb!{5.11} 
$c_1^*\o((f\o(Id_W\x c_2))\sptilde)\spcheck: \Bbb R \to 
C^\infty(\Bbb R,\Bbb R)$ is $C^\om$.
By corollary \nmb!{5.10} the associated mapping
$(c_1^*\o((f\o(Id_W\x c_2))\sptilde)\spcheck)\sptilde = 
c_2^*\o \check f \o c_1: \Bbb R \to  C^\om(\Bbb R,\Bbb R)$
is $C^\infty$. So for both families  describing the dual of
$C^\om(U,\Bbb R)$ we have shown that the composite with $\check
f\o c_1$ is $C^\infty$, so $\check f\o c_1$ is $C^\infty$.

If $c_1$ is $C^\om$, then $f\o (c_1\x Id_U):\Bbb R\x U \to 
W \x U \to  \Bbb R$ is $C^\om$. By step 1 the
associated mapping $(f\o(c_1\x Id_U))\spcheck= \check f\o c_1:
\Bbb R\to  C^\om(U,\Bbb R)$ is $C^\om$.

\smallskip
{\bf Step 3.} The general case. 
$$\align &\qquad f:W\x U \to  F\text{ is }C^\om \\
&\Leftrightarrow\quad \la\o f:W \x U \to  \Bbb R \text{ is }C^\om\text{
for all }\la\in F'\\
&\Leftrightarrow\quad (\la\o f)\spcheck = \la_*\o \check f:
W\to C^\om(U,\Bbb R)\text{ is }C^\om,\text{ by step 2,}\\
&\Leftrightarrow\quad \check f:W \to  C^\om(U,F)\text{ is }C^\om.\qed
\endalign$$
\enddemo

\heading\totoc \nmb0{6}. Consequences of cartesian closedness
\endheading

Among all those dual pairings on a fixed vector space $E$ that generate 
the same real analytic structure there is a finest one, namely that 
having as dual exactly the linear real analytic functionals, which are
exactly the bounded ones, by \nmb!{1.9}.
Recall that a dual pair $(E,E')$ is called convenient if and only if
it is complete and $E'$ consists exactly of the bounded linear 
functionals.

\proclaim{\nmb.{6.1}. Theorem}
The category of real analytic mappings between complete dual
pairs is equivalent 
to that of real analytic mappings between convenient dual pairs. Hence the
later category is also cartesian closed.
\endproclaim
\demo{Proof} The second category is a full subcategory of the first. A 
functor in the other direction is given by associating to every dual pair
$(E,E')$ the dual pair $(E,E^b)$, where 
$$\align E^b : 
&= \{\la:E \to  \Bbb R: \la\text{ is linear and bounded} \}\\
&= \{\la\in C^\om(E,\Bbb R): \la \text{ is linear}\}\\
&= \{\la:E \to  \Bbb R:
\la\text{ is linear, }
  \la\o c\in C^\om(\Bbb R,\Bbb R)\text{ for all }c\in C^\om(\Bbb R,E)\}.
\endalign$$
Functoriality follows since the real analytic mappings form a category. One composite
of this functor with the inclusion functor is the identity, and the other
is naturally isomorphic to the identity, since $E^b$ and $E'$ 
generate the same bornology and hence the same $C^\om$- and 
$C^\infty$-curves, by \nmb!{1.10}.
\qed
\enddemo

\subheading{Convention} All spaces are from now on assumed to be 
con\-ve\-nient and all function spaces will be considered with their 
natural bornological topology.

\proclaim{\nmb.{6.2}. Corollary (Canonical mappings are real analytic)}
The following mappings are $C^\om$:
\roster
\item $\ev:C^\om(U,F)\x U \to  F$, $(f,x)\mapsto f(x)$,
\item $\operatorname{ins}:E \to  C^\om(F,E\x F)$, $x\mapsto (y\mapsto (x,y))$,
\item $(\quad)^\wedge :C^\om(U,C^\om(V,G)) \to  C^\om(U\x V,G)$,
\item $(\quad)\spcheck:C^\om(U\x V,G) \to  C^\om(U,C^\om(V,G))$,
\item $\operatorname{comp}:C^\om(F,G)\x C^\om(U,F) \to  
    C^\om(U,G)$, $(f,g)\mapsto f\o g$,
\item $C^\om(\quad,\quad):C^\om(E_2,E_1)\x C^\om(F_1,F_2) \to $\newline
      $\to C^\om(C^\om(E_1,F_1),C^\om(E_2,F_2))$, 
      $(f,g)\mapsto (h\mapsto g\o h\o f)$.
\endroster
\endproclaim
\demo{Proof}
\therosteritem 1. The mapping associated to {\sl ev} via cartesian 
closedness is the identity on $C^\om(U,F)$, which is $C^\om$, 
thus {\sl ev} is also $C^\om$.

\therosteritem 2. The mapping associated to {\sl ins} via cartesian
closedness is the identity on $E\x F$, hence {\sl ins} is $C^\om$.

\therosteritem 3. The mapping associated via cartesian closedness 
is $(f;x,y)\mapsto f(x)(y)$, which is the $C^\om$-mapping $\ev\o(\ev\x 
id)$.

\therosteritem 4. The mapping associated by applying cartesian 
closedness twice is $(f;x;y)\mapsto f(x,y)$, which is just a 
$C^\om$ evaluation mapping.

\therosteritem 5. The mapping associated to {\sl comp}
via cartesian closedness is just $(f,g;x)\mapsto f(g(x))$, which 
is the $C^\om$-mapping $\ev\o (id \x \ev)$.

\therosteritem 6. The mapping associated by applying cartesian 
closed twice is $(f,g;h,x)\mapsto g(h(f(x)))$, which is the 
$C^\om$-mapping $\ev\o(id\x \ev)\o(id\x id\x \ev)$.
\qed
\enddemo

\proclaim{\nmb.{6.3}. Lemma (Canonical isomorphisms)}
One has the following natural isomorphisms:
\roster
\item $C^\om(W_1,C^\om(W_2,F))
    \cong C^\om(W_2,C^\om(W_1,F))$,
\item $C^\om(W_1,C^\infty(W_2,F))\cong C^\infty(W_2,C^\om(W_1,F))$.
\item $C^\om(W_1,L(E,F))\cong L(E,C^\om(W_1,F))$.
\item $C^\om(W_1,\ell^\infty(X,F))\cong \ell^\infty(X,C^\om(W_1,F))$.
\item $C^\om(W_1,\Cal Lip^k(X,F))\cong \Cal Lip^k(X,C^\om(W_1,F))$.
\endroster
In \therosteritem 4 $X$ is a $\ell^\infty$-space, i.e. a set together 
with a bornology induced by a family of real valued functions on $X$, cf. 
\cit!{5}, 1.2.4.
In \therosteritem 5 $X$ is a $\Cal Lip^k$-space, cf. 
\cit!{5}, 1.4.1.
The spaces $\ell^\infty(X,F)$ and $\Cal 
Lip^k(W,F)$ are defined in 
\cit!{5}, 3.6.1 and 4.4.1. 
\endproclaim
\demo{Proof}
All isomorphisms, as well as their inverse mappings, are 
given by the flip of coordinates: $f\mapsto \tilde f$, where 
$\tilde f(x)(y):=f(y)(x)$. Furthermore all occurring 
function spaces are convenient and satisfy the uniform $\Cal 
S$-boundedness theorem, where $\Cal S$ is the set of point 
evaluations, by \nmb!{5.5}, \nmb!{5.8}, \nmb!{4.6}, and 
by \cit!{5}, 3.6.1, 4.4.2, 3.6.6, and 4.4.7.

That $\tilde f$ has values in the corresponding spaces follows 
from the equation $\tilde f(x)=ev_x\o f$. 
One only has to check that $\tilde f$ itself 
is of the corresponding class, since it follows that $f\mapsto \tilde 
f$ is bounded. This is a consequence of the uniform boundedness 
principle, since 
$$(\ev_x\o \tilde{(\quad)})(f) = \ev_x(\tilde f)= 
\tilde f(x) = \ev_x\o f = (\ev_x)_*(f).$$ 

That $\tilde f$ is of the appropriate class in \therosteritem 1 and 
\therosteritem 2 follows by composing with 
$c_1 \in C^{\be_1}(\Bbb R,W_1)$ and 
$C^{\be_2}(\la,c_2):C^{\al_2}(W_2,F) \to  C^{\be_2}(\Bbb R,\Bbb R)$ for 
all $\la\in F'$ and $c_2\in C^{\be_2}(\Bbb R,W_2)$, where $\be_k$ and 
$\al_k$ are in $\{\infty,\om\}$ and $\be_k \leq \al_k$ for 
$k\in\{1,2\}$.
Then 
$C^{\be_2}(\la,c_2)\o \tilde f \o c_1 =
(C^{\be_1}(\la,c_1)\o f \o c_2)\sptilde
: \Bbb R \to  C^{\be_2}(\Bbb R,\Bbb R)$
is $C^{\be_1}$ by \nmb!{5.1} and \nmb!{5.10}, since 
$C^{\be_1}(\la,c_1)\o f \o c_2
: \Bbb R \to  W_2 \to  C^{\al_1}(W_1,F) \to  C^{\be_1}(\Bbb R,\Bbb R)$
is $C^{\be_2}$.

That $\tilde f$ is of the appropriate class in \therosteritem 3 follows,
since $L(E,F)$ is the $c^\infty$-closed subspace of $C^\om(E,F)$ formed
by the linear $C^\om$-mappings.

That $\tilde f$ is of the appropriate class in \therosteritem 4 follows from
\therosteritem 3, using the free convenient vector space $\ell^1(X)$ over 
the $\ell^\infty$-space $X$, see \cit!{5}, 5.1.24, 
satisfying $\ell^\infty(X,F)\cong L(\ell^1(X),F)$.

That $\tilde f$ is of the appropriate class in \therosteritem 5
follows from 
\therosteritem 3, using the free convenient vector space $\la^k(X)$ over 
the $\Cal Lip^k$-space $X$, see \cit!{5}, 5.1.3,
satisfying $\Cal Lip^k(X,F)\cong L(\la^k(X),F)$.
\qed
\enddemo

\subheading{Definition} A $C^{\infty,\om}$-mapping $f:U\x V \to  F$ is a 
mapping for which 
$$\check f\in C^\infty(U,C^\om(V,F)).$$

\proclaim{\nmb.{6.4}. Theorem (Composition of $C^{\infty,\om}$-mappings)}
Let $f:U\x V\to F$ and $g:U_1\x V_1\to V$ be $C^{\infty,\om}$, and 
$h:U_1\to U$ be $C^\infty$. Then 
$$f\o(h\o pr_1,g):U_1\x V_1\to F,\qquad
     (x,y)\mapsto f(h(x),g(x,y))$$
is $C^{\infty,\om}$.
\endproclaim
\demo{Proof}
We have to show that the mapping $x\mapsto (y\mapsto f(h(x),g(x,y))$, 
$U_1\to C^\om(V_1,F)$ is $C^\om$. It is well-defined, since $f$ and $g$ are 
$C^\om$ in the second variable. In order to show that it is 
$C^\om$ we compose with $\la_*:C^\om(V_1,F)\to C^\om(V_1,\Bbb R)$, 
where $\la\in F'$ is arbitrary. Thus it is enough to consider the 
case $F=\Bbb R$. Furthermore, we compose with $c^*:C^\om(V_1,\Bbb 
R)\to C^\al(\Bbb R,\Bbb R)$, where $c\in C^\al(\Bbb R,V_1)$ is 
arbitrary for $\al$ equal to $\om$ and $\infty$.

In case $\al=\infty$ the composite with  
$c^*$ is $C^\infty$, since the associated mapping $U_1\x \Bbb R\to \Bbb R$ is 
$f\o(h\o pr_1, g\o(id\x c))$ which is $C^\infty$.

Now the case $\al=\om$. Let $I\subseteq\Bbb R$ be an arbitrary 
open bounded interval. Then $c^*\o \check g:U_1\to C^\om(\Bbb R,G)$ 
is $C^\infty$, where $G$ is the convenient vector space 
containing $V$ as an $c^\infty$-open subset, and has values in 
the open set $\{\ga:\ga(\bar I)\subseteq V\}\subseteq C^\om(\Bbb 
R,G)$. Thus the composite with $c^*$, 
$\operatorname{comp}\o(\check f\o h,c^*\o \check g)$ is 
$C^\infty$, since
$\check f\o h:U_1\to U\to C^\om(V,F)$ is $C^\infty$,
$c^*\o \check g:U_1\to C^\om(\Bbb R,G)$ 
is $C^\infty$ and
$\operatorname{comp}:C^\om(V,F)\x\{\ga\in C^\om(\Bbb R,G):\ga(\bar 
I)\subseteq V\}\to C^\om(I,\Bbb R)$ is $C^\om$, because it is associated to
$\ev\o(id\x \ev):C^\om(V,F)\x\{\ga\in C^\om(\Bbb R,G):\ga(\bar 
I)\subseteq V\}\x I\to \Bbb R$. That $\ev:\{\ga\in C^\om(\Bbb R,G):\ga(\bar 
I)\subseteq V\}\x I\to \Bbb R$ is $C^\om$ follows, since
the associated mapping is the restriction mapping $C^\om(\Bbb R,G)\to 
C^\om(I,G)$.
\qed\enddemo

\proclaim{\nmb.{6.5} Corollary}
Let $f:U\to  F$ be $C^\om$ and $g:U_1\x V_1\to  U$ be $C^{\infty,\om}$,
then 
$$f\o g:U_1\x V_1\to  F$$
is $C^{\infty,\om}$.

Let $f:U\x V\to  F$ be $C^{\infty,\om}$ and $h:U_1\to  U$ be 
$C^\infty$, then $f\o(h\x id):U_1\x V\to  F$ is $C^{\infty,\om}$.\qed
\endproclaim

The second part is a generalization of theorem \nmb!{5.11}.

\proclaim{\nmb.{6.6}. Corollary}
Let $f:E\supseteq U\to F$ be $C^\om$, let $I\subseteq \Bbb R$ be 
open and bounded, and $\al$ be $\om$ or $\infty$. Then 
$f_*:C^\al(\Bbb R,E)\supseteq 
\{c:c(\bar I)\subseteq U\} \to C^\al(I,F)$
is $C^\om$.
\endproclaim
\demo{Proof} Obviously $f_*(c):=f\o c\in C^\al(I,F)$ is 
well-defined for all $c\in C^\al(\Bbb R,E)$ satisfying 
$c(\bar I)\subseteq U$.

Furthermore the composite of $f_*$ with any $C^\be$-curve 
$\ga:\Bbb R\to \{c:c(\bar I)\subseteq U\}\subseteq C^\al(\Bbb R,E)$ 
is a $C^\be$-curve in $C^\al(I,F)$ for $\be$ equal to $\om$ or 
$\infty$. For $\be=\al$ this follows from cartesian closedness
of the $C^\al$-maps. For $\al\not = \be$ this follows from 
\nmb!{6.5}.

Finally $\{c:c(\bar I)\subseteq U\}\subseteq 
C^\al(\Bbb R,E)$ is $c^\infty$-open, since it is open for the 
topology of uniform convergence on compact sets which is coarser
than the bornological and hence than the $c^\infty$-topology on 
$C^\al(\Bbb R,E)$. Here is the only place where we make use of the 
boundedness of $I$.
\qed\enddemo

\proclaim{\nmb.{6.7}. Lemma (Free convenient vector space)}
Let $U\subseteq E$ be $c^\infty$-open in a convenient vector space 
$E$. There exists a free convenient vector spaces $Free(U)$ over 
$U$, i.e. for every convenient vector 
space $F$, one has a natural isomorphism $C^\om(U,F)\cong 
L(Free(U),F)$
\endproclaim
\demo{Proof}
Consider the Mackey closure $Free(U)$ of the linear subspace of 
$C^\om(U,\Bbb R)'$ generated by the set $\{ev_x:x\in U\}$. Let 
$\io:U\to Free(U)\subseteq C^\om(U,\Bbb R)'$ be the mapping given 
by $x\mapsto \ev_x$. This mapping is $C^\om$, since $ev_f\o\io=f$ 
for every $f\in C^\om(U,\Bbb R)$.

Obviously every real analytic mapping $f:U\to F$ extends to the linear bounded
mapping $\tilde f:C^\om(U,\Bbb R)' \to  F''$,
$\la \mapsto (l \mapsto \la(l\o f))$. Since $\tilde f$
coincides on the generators $ev_x$ with 
$f$, it maps the Mackey closure $Free(U)$ into the Mackey 
closure of $F\supseteq f(U)$. Since $F$ is complete this is again 
$F$. Uniqueness of $\tilde f$ follows, since every
linear real analytic mapping is bounded, hence it is determined by its values on the 
subset $\{\ev_x:x\in U\}$ that spans the linear subspace having as 
Mackey closure $Free(U)$.
\qed
\enddemo

\proclaim{\nmb.{6.8}. Lemma (Derivatives)}
The derivative $d$, where 
$$df(x)(v):=\frac d{dt}\mid_{t=0}f(x+tv),$$
is bounded and linear $d:C^\om(U,F) \to  C^\om(U,L(E,F))$.
\endproclaim
\demo{Proof}
The differential $df(x)(v)$ makes
sense and is linear in $v$, because every real analytic mapping $f$ is
smooth. So it remains to show that 
$(f,x,v)\mapsto df(x)(v)$ is real analytic. So let $f$, $x$, and $v$ depend 
real analytically (resp. smoothly) on a real parameter $s$. 
Since $(t,s)\mapsto x(s)+tv(s)$ is real analytic (resp. smooth) into $U\subseteq E$, 
the mapping $r\mapsto ((t,s)\mapsto f(r)(x(s)+tv(s))$ is real analytic into
$C^\om(\Bbb R^2,F)$ (resp. $C^\infty(\Bbb R^2,F$). Composing with
$$\frac d{dt}\mid_{t=0}:C^\om(\Bbb R^2,F) \to  C^\om(\Bbb R,F)\qquad
(\text{resp. } :C^\infty(\Bbb R^2,F) \to  C^\infty(\Bbb R,F))$$ 
shows that
$r\mapsto (s\mapsto d(f(r))(x(s))(v(s)))$, $\Bbb R \to  C^\om(\Bbb R,F)$
is real analytic. Considering the associated mapping on $\Bbb R^2$ composed
with the diagonal map shows that $(f,x,v)\mapsto df(x)(v)$ is real analytic. 
\qed
\enddemo

The following examples as well as several others can be found in 
\cit!{5}, 5.3.6.

\proclaim{\nmb.{6.9} Example} Let
$T:C^\infty(\Bbb R,\Bbb R) \to  C^\infty(\Bbb R,\Bbb R)$ be given by 
$T(f)=f'$. Then
the continuous linear differential equation
$x'(t)=T(x(t))$ with initial value $x(0)=x_0$
has a unique smooth solution $x(t)(s)=x_0(t+s)$ which is 
however not real analytic.
\endproclaim
\demo{Proof}
A smooth curve $x:\Bbb R \to  C^\infty(\Bbb R,\Bbb R)$ is a solution of
the differential equation $x'(t)=T(x(t))$ if and only if 
$$\frac{\partial}{\partial t}\hat x(t,s) = 
\frac{\partial}{\partial s}\hat x(t,s).$$
Hence we have
$$\frac{d}{dt}\hat x(t,r-t)=0,$$
i.e\. $\hat x(t,r-t)$ is constant and 
hence equal to $\hat x(0,r)=x_0(r)$. Thus $\hat x(t,s)=x_0(t+s)$.

Suppose $x:\Bbb R \to  C^\infty(\Bbb R,\Bbb R)$ were real analytic. Then the 
composite with 
$$ev_0:C^\infty(\Bbb R,\Bbb R) \to  \Bbb R$$
were a real 
analytic function. But this composite is just $x_0=ev_0\o x$, which is 
not in general real analytic.
\qed\enddemo

\proclaim{\nmb.{6.10} Example}
Let $E$ be either $C^\infty(\Bbb R,\Bbb R)$ or 
$C^\om(\Bbb R,\Bbb R)$. Then the mapping $exp_*:E \to  E$ is $C^\om$, has 
invertible derivative at every point, but the image does not contain an open 
neighborhood of $exp_*(0)$.
\endproclaim
\demo{Proof}
That $exp_*$ is $C^\om$ was shown in \nmb!{6.6}. Its derivative is 
given by 
$$(exp_*)'(f)(g):t\mapsto g(t)e^{f(t)}$$
and hence is invertible 
with $g\mapsto (t\mapsto g(t)e^{-f(t)})$ as inverse mapping.
Now consider the real analytic curve $c:\Bbb R \to  E$ given by 
$c(t)(s)=1-(ts)^2$. One has $c(0)=1=exp_*(0)$, but $c(t)$ is not in the 
image of $exp_*$ for any $t\ne 0$, since $c(t)(\frac1t)=0$ but 
$exp_*(g)(t)=e^{g(t)}>0$ for all $g$ and $t$.
\qed\enddemo

\heading\totoc \nmb0{7}. Spaces of sections of vector bundles \endheading

\subheading{\nmb.{7.1}. Vector bundles} Let $(E,p,M)$ be a real analytic finite
dimensional vector bundle over a real analytic manifold $M$,
where $E$ is their total space and $p:E\to M$ is the projection.
So there is an open cover $(U_\al)_\al$ of $M$ and vector bundle charts
$\ps_\al$ satisfying
$$\CD
E| U_\al := p\i(U_\al)  @>\ps_\al>>  U_\al\x V \\
@VpVV                               @VV{pr_1}V \\
U_\al                          @=       U_\al.
\endCD$$
Here $V$ is a fixed finite dimensional real vector space, called
the standard fiber.
We have $(\ps_\al\o\ps_\be\i)(x,v) = (x,\ps_{\al\be}(x)v)$ for transition
functions $\ps_{\al\be}:U_{\al\be} = U_\al\cap U_\be \to GL(V)$,
which are real analytic.

If we extend the transition functions $\ps_{\al\be}$ to 
$\widetilde\ps_{\al\be}:\widetilde U_{\al\be} \to GL(V_{\Bbb C}) =
GL(V)_{\Bbb C}$, we see that there is a holomorphic vector
bundle $(E_{\Bbb C},p_{\Bbb C},M_{\Bbb C})$ over a complex (even
Stein) manifold $M_{\Bbb C}$ such that $E$ is isomorphic to a real part
of $E_{\Bbb C}| M$, compare \nmb!{3.10}.
The germ of it along $M$ is unique.

Real analytic sections $s:M\to E$ coincide with certain germs
along $M$ of holomorphic sections $W \to E_{\Bbb C}$ for open
neighborhoods $W$ of $M$ in $M_{\Bbb C}$.

\subheading{\nmb.{7.2}. Spaces of sections} For a holomorphic vector
bundle $(F,q,N)$ over a complex manifold $N$ we denote by $\Cal
H(F)$ the vector space of all holomorphic sections $s:N\to F$,
equipped with the compact open topology, a nuclear Fr\'echet
topology, since it is initial with respect to the cone
$$\Cal H(F)\to \Cal H(F|U_\al) @>{(pr_2\o\ps_\al)_*}>>
\Cal H(U_\al,\Bbb C^k) = \Cal H(U_\al,\Bbb C)^k,$$
of mappings into nuclear spaces, see \nmb!{3.2}.

For a subset
$A\subseteq N$ let $\Cal H(F| A)$ be the space of
germs along $A$ of holomorphic sections  $W\to F| W$ for open
sets $W$ in $N$ containing
$A$. We equip $\Cal H(F| A)$ with the locally convex topology
induced by the inductive cone 
$\Cal H(F| W)\to \Cal H(F| A)$ for all $W$.
This is Hausdorff since jet prolongations at points in $A$
separate germs. 

For a real analytic vector bundle $(E,p,M)$
let $C^\omega(E)$ be the space of real analytic sections
$s:M\to E$.
Furthermore let $C^\om(E| A)$ denote the space of
germs at a subset $A\subseteq M$ of real analytic sections
defined near $A$.
The complexification of this real vector space
is the complex vector space $\Cal H(E_{\Bbb C}| A)$,
because germs of real analytic sections $s:A\to E$ extend
uniquely to germs
along $A$ of holomorphic sections $W \to E_{\Bbb C}$ for open
sets $W$ in $M_{\Bbb C}$ containing $A$, compare \nmb!{3.11}. 

We topologize $C^\om(E| A)$ 
as subspace of $\Cal H(E_{\Bbb C}| A)$.

For a smooth
vector bundle $(E,p,M)$ let $C^\infty(E)$ denote the nuclear Fr\'echet
space of all smooth sections with the topology of uniform
convergence on compact subsets, in all derivatives separately, see
\cit!{18} and \cit!{5}, 4.6.

\proclaim{\nmb.{7.3}. Theorem (Structure on spaces of germs of sections)} 
If $(E,p,M)$ is a real analytic vector bundle and $A$ a closed
subset of $M$,
then the space $C^\om(E| A)$ is convenient.
Its bornology is generated by the cone 
$$C^\om(E| A) @>{(\ps_\al)_*}>> 
C^\om(U_\al\cap A\subseteq U_\al,\Bbb R)^k,$$
where $(U_\al,\ps_\al)_\al$ is an arbitrary real analytic vector
bundle atlas of $E$. If $A$ is compact, the space $C^\om(E|A)$
is nuclear.
\endproclaim
The corresponding statement for smooth sections is also true, see
\cit!{5}, 4.6.23. 
\demo{Proof} We show the corresponding result for
holomorphic germs. By taking real parts the theorem then
follows.  

So let $(F,q,N)$ be a holomorphic vector
bundle and let $A$ be a closed subset of $N$.
Then $\Cal H(F| A)$ is a bornological
locally convex space, since it is an inductive limit of the
Fr\'echet spaces $\Cal H(F| W)$ for open sets $W$
containing $A$. If $A$ is compact, $\Cal H(F|A)$ is nuclear as
countable inductive limit. 

Let $(U_\al,\ps_\al)_\al$ be a holomorphic
vector bundle atlas for $F$. 

Then we consider the cone 
$$\Cal H(F| A) @>{(\ps_\al)_*}>> 
\Cal H(U_\al\cap A\subseteq U_\al,\Bbb C^k) = 
\Cal H(U_\al\cap A\subseteq U_\al,\Bbb C)^k.$$
Obviously each mapping is continuous, so the cone induces a
bornology which is coarser than the given one, and which is
complete by \nmb!{3.13}. 

It remains to show that every subset 
$\Cal B\subseteq \Cal H(F| A)$, such that 
$(u_\al)_*(\Cal B)$ is bounded in every  
$\Cal H(U_\al\cap A\subseteq U_\al,\Bbb C)^k$, is bounded in $\Cal
H(F| W)$ for some open neighborhood $W$ of $A$ in $N$.

Since all restriction mappings to
smaller subsets are continuous, it suffices to show the
assertions of the theorem for some refinement of the atlas
$(U_\al)$. Let us pass first to a relatively compact refinement.
By topological dimension theory there is a further refinement 
such that any $\dim_{\Bbb R}N+2$ different sets have empty
intersection. We call the resulting atlas again $(U_\al)$.
Let $(K_\al)$ be a cover of $N$ consisting of compact subsets
$K_\al\subseteq U_\al$ for all $\al$.

For any finite set $\Cal A$ of indices
let us consider now all non empty intersections
$U_{\Cal A}:=\bigcap_{\al\in\Cal A}U_\al$ and 
$K_{\Cal A}:=\bigcap_{\al\in\Cal A}K_\al$. 
Since by \nmb!{3.4} (or \nmb!{3.6}) the space
$\Cal H(A\cap K_{\Cal A}\subseteq U_{\Cal A},\Bbb C)$ 
is a regular inductive limit there are open sets 
$W_{\Cal A}\subseteq U_{\Cal A}$ containing $A\cap K_{\Cal A}$,
such that $\Cal B|(A\cap K_{\Cal A})$ (more precisely
$(\ps_{\Cal A})_*(\Cal B|(A\cap K_{\Cal A}))$ for some suitable
vector bundle chart mappings $\ps_{\Cal A}$) is contained and
bounded in $\Cal H(W_{\Cal A},\Bbb C)^k$.
By passing to smaller open sets we may assume that 
$W_{\Cal A_1}\subseteq W_{\Cal A_2}$ for $\Cal A_1\supseteq\Cal A_2$.
Now we define the subset 
$$W := \bigcup_{\Cal A}\widehat W_{\Cal A},\text{ where }
\widehat W_{\Cal A} := W_{\Cal A}\setminus\bigcup_{\al\notin\Cal A}K_\al$$
$W$ is open since $(K_\al)$ is a locally finite cover.
For $x\in A$ let $\Cal A:=\{\al:x\in K_\al\}$, then $x\in \widehat W_{\Cal A}$.

Now we show that every germ $s\in \Cal B$ has a unique extension to $W$.
For every $\Cal A$ the germ of $s$ along $A\cap K_{\Cal A}$ has a unique
extension $s_{\Cal A}$ to a section over $W_{\Cal A}$ and for
$\Cal A_1\subseteq \Cal A_2$ we have $s_{\Cal A_1}|W_{\Cal A_2}
= s_{\Cal A_2}$. 
We define the extension $s_W$ by $s_W|\widehat
W_{\Cal A}=s_{\Cal A}|\widehat W_{\Cal A}$. This is well defined since
one may check that 
$\widehat W_{\Cal A_1}\cap\widehat W_{\Cal A_2}\subseteq 
\widehat W_{\Cal A_1\cap\Cal A_2}$.

$\Cal B$ is bounded in $\Cal H(F|W)$ if it is uniformly
bounded on each compact subset $K$ of $W$. This is true since each $K$ is
covered by finitely many $W_\al$ and $\Cal B|A\cap K_\al$ is bounded
in $\Cal H(W_\al,\Bbb C)$. 
\qed\enddemo

\subheading{\nmb.{7.4}} Let $f:E\to E'$ be a real analytic vector
bundle homomorphism, i.e. we have a commutative diagram
$$\CD
E  @>f>> E' \\
@VpVV    @VVp'V\\
M @>>\underline f> M' \endCD$$
of real analytic mappings such that $f$ is fiberwise linear. 

\proclaim{Lemma} If $f$ is fiberwise invertible, then
$f^*:C^\om(E')\to C^\om(E)$, given by 
$$(f^*s)(x) := (f_x)\i(s(\underline f(x)),$$ 
is continuous and linear.

If $\underline f = Id_M$ then the mapping $f_*:C^\om(E)\to
C^\om(E')$, given by $s\mapsto f\o s$, is continuous and linear.
\endproclaim
\demo{Proof} Extend $f$ to the complexification. Here for the
compact open topology on the corresponding spaces of holomorphic
sections the assertion is trivial.
\qed\enddemo

\subheading{\nmb.{7.5} Real analytic mappings are dense} Let
$(E,p,M)$ be a real analytic vector 
bundle. Then there is another real analytic vector bundle
$(E',p',M)$ such that the Whitney sum $E\oplus E'\to M$ is real
analytically isomorphic to a trivial bundle $M\x \Bbb R^k\to M$.
This is seen as follows:
By \cit!{7}, Theorem 3, there is a closed real
analytic embedding $i:E\to \Bbb R^k$ for some $k$. 
Now the fiber derivative along the zero section gives a
fiberwise linear and injective real analytic mapping $E\to
\Bbb R^k$ which induces a real analytic embedding $j$ of the vector
bundle $(E,p,M)$ into the trivial bundle $M\x \Bbb R^k\to M$.
The standard inner product on $\Bbb R^k$ gives rise to the real
analytic orthogonal complementary vector bundle $E':=E^\bot$ and
a real analytic Riemannian metric on the vector bundle $E$.

Hence an embedding of the real analytic vector bundle
into another one induces a linear embedding of the spaces of real
analytic sections onto a direct summand.

We remark that in this situation the orthogonal projection onto
the vertical bundle $VE$ within $T(M\x \Bbb R^k)$ gives rise to
a real analytic linear connection (covariant derivative)
$\nabla:C^\om(TM)\x C^\om(E) \to C^\om(E)$.
If $c:\Bbb R\to M$ is a smooth or real analytic curve in $M$, we
have the parallel transport
$Pt(c,t)v\in E_{c(t)}$ for all $v\in E_{c(0)}$ and $t\in \Bbb R$ 
which is smooth or real analytic, respectively, on $\Bbb R\x E_{c(0)}$.
It is given by the differential equation 
$\nabla_{\partial_t} Pt(c,t)v = 0$.

\proclaim{\nmb.{7.6} Corollary} If $\nabla$ is a real analytic linear
connection on a vector bundle $(E,p,M)$, then the following cone
generates the bornology on $C^\om(E)$.
$$C^\om(E) @>{Pt(c,\quad)^*}>> C^\al(\Bbb R,E_{c(0)}),$$ 
for all $c\in C^\al(\Bbb R,M)$ and $\al= \om,\infty$.
\endproclaim
\demo{Proof} The bornology induced by the cone is coarser that
the given one by \nmb!{7.4}. A still coarser bornology is induced
by all curves subordinated to some vector bundle atlas. Hence by
theorem \nmb!{7.3} it suffices to check for a trivial bundle, that
this bornology coincides with the given one.
So we assume that $E$ is trivial. For the constant parallel
transport the result follows from lemma \nmb!{5.3}. The change
to an arbitrary real analytic parallel transport can be
absorbed into a $C^\al$-isomorphism of each vector bundle $c^*E$
separately.
\qed\enddemo

\proclaim{\nmb.{7.7}. Lemma (Curves in spaces of sections)} \newline
1. For a real analytic vector bundle $(E,p,M)$ 
a curve $c:\Bbb R\to C^\om(E)$ is real analytic if
and only if the associated mapping
$\hat c:\Bbb R\x M \to E$ is real analytic.

The curve $c:\Bbb R\to C^\om(E)$ is smooth if and only
if $\hat c:\Bbb R\x M\to E$ satisfies
the following condition:
\roster
\item"" For each $n$ there is an open neighborhood $U_n$ of
$\Bbb R\x M$ in $\Bbb R\x M_{\Bbb C}$ and a (unique)
$C^n$-extension $\tilde c: U_n\to E_{\Bbb C}$ such that 
$\tilde c(t,\quad)$ is holomorphic for all $t\in \Bbb R$.
\endroster
2. For a smooth vector bundle $(E,p,M)$ a curve $c:\Bbb R\to
C^\infty(E)$ is smooth if and only if $\hat c:\Bbb R\x M\to E$
is smooth. 

The curve $c:\Bbb R\to C^\infty(E)$ is real analytic if and only
if $\hat c$ satisfies the following condition:
\roster
\item"" For each $n$ there is an open neighborhood $U_n$ of 
$\Bbb R\x M$ in $\Bbb C \x M$ and a (unique) $C^n$-extension $\tilde
c:U_n\to E_{\Bbb C}$ such that $\tilde c(\quad,x)$ is
holomorphic for all $x\in M$.
\endroster
\endproclaim
\demo{Proof} 1. By theorem \nmb!{7.3} we
may assume that $M$ is open in $\Bbb R^m$, and we may consider 
$C^\om(M,\Bbb R)$ instead of $C^\om(E)$. The statement on
real analyticity follows from cartesian closedness, \nmb!{5.12}.

To prove the statement on smoothness we note that $C^\om(M,\Bbb R)$
is the real part of $\Cal H(M\subseteq\Bbb C^m,\Bbb C)$ by
\nmb!{3.11}, which is a regular inductive limit of spaces
$\Cal H(W,\Bbb C)$ for open neighborhoods $W$ of $M$ in
$\Bbb C^m$ by \nmb!{3.6}. By \cit!{13}, Folgerung on p\. 114, 
$c$ is smooth if and only if for each $n$, locally in $\Bbb R$ it
factors to a $C^n$-curve into some $\Cal H(W,\Bbb C)$, which
sits continuously embedded in $C^\infty(W,\Bbb R^2)$. So the
associated mapping $\Bbb R\x M_{\Bbb C}\supseteq J\x W\to \Bbb C$ is $C^n$
and holomorphic in the second variables, and conversely.

2. By \nmb!{7.3} we
may assume that $M$ is open in $\Bbb R^m$, and we may consider 
$C^\infty(M,\Bbb R)$ instead of $C^\infty(E)$. The statement on
smoothness follows from cartesian closedness 
of smooth mappings, similarly as the $C^\om$-statement above.

To prove the statement on real analyticity we note that $C^\infty(M,\Bbb R)$
is the projective limit of the Banach spaces $C^n(M_i,\Bbb R)$, where
$M_i$ is a covering of $M$ by compact cubes. By lemma
\nmb!{1.11} the curve $c$ is real analytic if and only if it is
real analytic into each $C^n(M_i,\Bbb R)$, by \nmb!{1.6} and
\nmb!{1.7} it extends locally to a holomorphic curve 
$\Bbb C\to C^n(M_i,\Bbb C)$. Its associated mappings fit
together to the required $C^n$-extension $\tilde c$. 
\qed\enddemo

\proclaim{\nmb.{7.8}. Corollary} Let $(E,p,M)$ and $(E',p',M)$ be
real analytic vector bundles over a compact manifold $M$. Let
$W\subseteq E$ be an open subset such that $p(W) = M$, and let
$f:W\to E'$ be a fiber respecting real analytic (nonlinear) mapping.

Then $C^\infty(W):=\{s\in C^\infty(E): s(M)\subseteq W\}$ is
open and not empty in the convenient vector space $C^\infty(E)$.
The mapping $f_*:C^\infty(W)\to C^\infty(E')$ is real analytic
with derivative
$(d_vf)_*:C^\infty(W)\x C^\infty(E)\to C^\infty(E')$, where the vertical
derivative $d_vf:W\x_ME\to E'$ is given by 
$$d_vf(u,w):= \left.\frac d{dt}\right|_0 f(u+tw).$$

Then $C^\om(W):= \{s\in C^\om(E): s(M)\subseteq W\}$ is open and
not empty in the convenient vector space $C^\om(E)$ and the
mapping $f_*:C^\om(W)\to C^\om(E')$ is real analytic with derivative
$(d_vf)_*:C^\om(W)\x C^\om(E)\to C^\om(E')$.
\endproclaim

\demo{Proof} The set
$C^\infty(W)$ is open in $C^\infty(E)$ since it is open in the
compact-open topology. Then $C^\om(W)$
is open in $C^\om(E)$ since $C^\om(E)\to C^\infty(E)$ is
continuous by \nmb!{3.12} and \nmb!{7.3}. 

Now we prove the statement for $C^\om(W)$, the proof for
$C^\infty(W)$ is then similar.

We check that $f_*$ maps $C^\om$-curves to $C^\om$-curves and
maps smooth curves to smooth curves.

If $c:\Bbb R\to C^\om(W)\subseteq C^\om(E)$ is $C^\om$, then 
$\hat c:\Bbb R\x M\to E$ is $C^\om$ by lemma \nmb!{7.7}. So
$(f_*\o c)\sphat=f\o\hat c:\Bbb R\x M\to E'$ is also $C^\om$,
hence $f_*\o c:\Bbb R \to C^\om(E')$ is $C^\om$.

If $c:\Bbb R\to C^\om(W)\subseteq C^\om(E)$ is smooth, for each
$n$ there is an open neighborhood $U_n$ of $\Bbb R\x M$ in 
$\Bbb R\x M_{\Bbb C}$ and a $C^n$-extension $\tilde c:U_n\to
E_{\Bbb C}$ of $\hat c$ such that $\tilde c(t,\quad)$ is
holomorphic. The mapping $f:W\to E'$ has also a holomorphic
extension $\tilde f:E_{\Bbb C}\supseteq W_{\Bbb C}\to E'_{\Bbb C}$.
Then $\tilde f\o \tilde c$ is an extension of $(f_*\o c)\sphat$
satisfying the condition in lemma \nmb!{7.7}, so 
$f_*\o c:\Bbb R\to C^\om(E')$ is smooth.
\qed\enddemo

\heading\totoc \nmb0{8}. Manifolds of analytic mappings \endheading

\subheading{\nmb.{8.1}. Infinite dimensional real analytic manifolds}
A chart $(U,u)$ on a set $\Cal M$ is a bijection $u:U\to
u(U)\subseteq E_U$ from a subset $U\subseteq \Cal M$ onto a
$c^\infty$-open subset of a convenient vector space $E_U$. Two
such charts are called $C^\om$-compatible, if the chart change
mapping $u\o v\i:v(U\cap V) \to u(U\cap V)$ is a  $C^\om$-diffeomorphism
between $c^\infty$-open subsets of convenient vector spaces.
A $C^\om$-atlas on $\Cal M$ is a set of pairwise
$C^\om$-compatible charts on $\Cal M$ which cover $\Cal M$. 
Two such atlases are equivalent if their union is again a
$C^\om$-atlas. A $C^\om$-structure on $\Cal M$ is an equivalence
class of $C^\om$-atlases. A $C^\om$-manifold $\Cal M$ is a set
together with a $C^\om$-structure on it.

The natural topology on $\Cal M$ is the 
identification topology, where a subset
$W\subseteq \Cal M$ is open if and only if $u(U\cap W)$ is
$c^\infty$-open in $E_U$ for all charts in a $C^\om$-atlas
belonging to the structure. 
 In the finite dimensional treatment
of manifolds one requires that this topology has some
properties: Hausdorff, separable or metrizable or paracompact.

In infinite dimensions it is not yet clear what the most
sensible requirements are. Hausdorff does not imply regular. If
$\Cal M$ is Hausdorff and regular, and if all modeling vector
spaces admit smooth bump functions, any locally finite open
cover of $\Cal M$ admits a subordinated smooth partition of unity. 

Mappings between $C^\om$-manifolds are called $C^\infty$ or
$C^\om$ if they are continuous and their chart representations
are smooth or real analytic, respectively. 

The final topology
with respect to all smooth 
curves coincides with the identification topology on a
$C^\om$-manifold. So the following two statements hold:

A mapping $f:\Cal M\to \Cal N$ between
$C^\om$-manifolds is $C^\infty$ if $f\o c$ is $C^\infty$ for
each $C^\infty$-curve in $\Cal M$.

$f$ is $C^\om$ if it is $C^\infty$ and $f\o c$ is $C^\om$ for
each $C^\om$-curve in $\Cal M$. 

The tangent bundle $T\Cal M\to \Cal M$ of a $C^\om$-manifold
$\Cal M$ is the vector bundle glued from the sets $u(U)\x E_U$ via the
transition functions 
$(x,y)\mapsto ((u\o v\i)(x),d(u\o v\i)(x)y)$ for all charts
$(U,u)$ and $(V,v)$ in a $C^\om$-atlas of $\Cal M$.

\proclaim{\nmb.{8.2}. Theorem (Manifold structure of
$C^\om(M,N)$)} Let $M$ and $N$ be real analytic 
manifolds, let $M$ be compact. Then the space $C^\om(M,N)$ of
all real analytic mappings from $M$ to $N$ is a real analytic
manifold, modeled on spaces $C^\om(f^*TN)$ of real analytic
sections of pullback bundles along $f:M\to N$ over $M$.
\endproclaim
\demo{Proof} Choose a real analytic Riemannian metric on $N$.
See \nmb!{7.5} for a sketch how to find one. Let $\exp: TN\supseteq
U\to N$ be the real analytic exponential mapping of this Riemannian metric,
defined on a suitable open neighborhood of the zero section. We
may assume that $U$ is chosen in such a way that
$(\pi_N,\exp):U\to N\x N$ is a real analytic diffeomorphism onto
an open neighborhood $V$ of the diagonal.

For $f\in C^\om(M,N)$ we consider the pullback vector bundle
$$\CD M\x_NTN @= f^*TN     @>{\pi_N^*f}>>           TN \\
@.               @V{f^*\pi_N}VV                @VV{\pi_N}V\\
              @. M         @>f>>              N.  \endCD$$
Then $C^\om(f^*TN)$ is canonically isomorphic to the space
$$C^\om_f(M,TN):= \{h\in C^\om(M,TN):\pi_N\o h=f\}$$ 
via $s\mapsto (\pi_N^*f)\o s$ and $(Id_M,h)\gets h$.

Now let 
$$U_f :=\{g\in C^\om(M,N):(f(x),g(x))\in V\text{ for all }x\in M\}$$ 
and let $u_f:U_f\to C^\om(f^*TN)$ be given by 
$$u_f(g)(x) = (x,\exp_{f(x)}\i(g(x))) = (x,((\pi_N,\exp)\i\o(f,g))(x)).$$
Then $u_f$ is a bijective mapping from $U_f$ onto 
$\{s\in C^\om(f^*TN): s(M)\subseteq f^*U\}$, whose inverse is
given by $u_f\i(s) = \exp\o(\pi_N^*f)\o s$, where we view 
$U \to N$ as fiber bundle. Since $M$ is compact, $u_f(U_f)$ is open in
$C^\om(f^*TN)$ for the compact open topology, thus also for the
finer topology described in \nmb!{7.2}.

Now we consider the atlas $(U_f,u_f)_{f\in C^\om(M,N)}$ for
$C^\om(M,N)$. Its chart change mappings are given 
for $s\in u_g(U_f\cap U_g)\subseteq C^\om(g^*TN)$ by 
$$\align (u_f\o u_g\i)(s) &= (Id_M,(\pi_N,\exp)\i\o(f,\exp\o(\pi_N^*g)\o s)) \\
&= (\tau_f\i\o\tau_g)_*(s),\endalign$$ 
where $\tau_g(x,Y_{g(x)}) := (x,\exp_{g(x)}(Y_{g(x)})))$
is a real analytic diffeomorphism 
$$\tau_g:g^*TN \supseteq g^*U \to (g\x Id_N)\i(V)\subseteq M\x N$$
which is fiber respecting over $M$. Thus by \nmb!{7.8} the chart change
$u_f\o u_g\i = (\tau_f\i\o \tau_g)_*$ is defined on an open
subset and real analytic.

Finally we put the identification topology from this atlas onto
the space $C^\om(M,N)$, which is obviously finer than the
compact open topology and thus Hausdorff.

The equation $u_f\o u_g\i = (\tau_f\i\o \tau_g)_*$ shows that
the real analytic structure does not depend on the choice of the
real analytic Riemannian metric on $N$.
\qed\enddemo

\demo{Remark} If $N$ is a finite dimensional vector space, then
the structure of a real analytic manifold on $C^\om(M,N)$
described here coincides with that of the space 
$C^\om(M\x N\to M)$ of sections discussed in
\nmb!{7.2}, because the exponential mapping of any euclidean
structure is the affine structure of $N$.
\enddemo

\proclaim{\nmb.{8.3}. Theorem ($C^\om$-manifold structure on
$C^\infty(M,N)$)} Let $M$ and $N$ be real 
analytic manifolds, with $M$ compact. Then the smooth manifold
$C^\infty(M,N)$ with the structure from \cit!{18}, 10.4 
is even a real analytic manifold. 
\endproclaim
\demo{Proof} For a fixed real analytic exponential mapping on $N$
the charts $(U_f,u_f)$ (from \nmb!{8.2} with $C^\om$ replaced
by $C^\infty$, which agrees with those from \cit!{18}, 10.4, 
see also \cit!{5}, 4.7) for $f\in C^\om(M,N)$ form a smooth atlas
for $C^\infty(M,N)$, since $C^\om(M,N)$ is dense in
$C^\infty(M,N)$ by \cit!{7}, Proposition 8. The chart
changings $u_f\o u_g\i = (\tau_f\i\o \tau_g)_*$ are real
analytic by \nmb!{7.8}. 
\qed\enddemo

\subheading{\nmb.{8.4} Remark} If $M$ is not compact,
$C^\om(M,N)$ is dense in $C^\infty(M,N)$ for the
Whitney-$C^\infty$-topology by \cit!{7}, Proposition 8.
This is not the case for the $\Cal D$-topology from \cit!{18}, in which
$C^\infty(M,N)$ is a smooth manifold. The charts $U_f$ for $f\in
C^\om(M,N)$ do not cover $C^\infty(M,N)$.

\proclaim{\nmb.{8.5}. Theorem} Let $M$
and $N$ be real analytic manifolds, where $M$ is compact, the
two infinite dimensional real analytic 
vector bundles $TC^\om(M,N)$ and $C^\om(M,TN)$ over $C^\om(M,N)$
are canonically isomorphic. The same assertion is true for
$C^\infty(M,N)$. 
\endproclaim
\demo{Proof} Let us fix an exponential mapping $\exp$ on $N$.
It gives rise to the canonical atlas $(U_f,u_f)$ for
$C^\om(M,N)$ from \nmb!{8.2}. 
$TC^\om(M,N)$ is defined as the vector bundle glued
from the transition functions $(r,s)\mapsto
(u_f(u_g\i(r)),d(u_f\o u_g\i)(r)s)$.
Then $T(\exp)$ composed with the canonical flip on $T^2N$ is an
exponential mapping for $TN$, which gives rise to the canonical
atlas $(U_{0\o f},u_{0\o f})$ for $C^\om(M,TN)$, where $0$ is
the zero section of $TN$. Via some
canonical identifications the 
two sets of transition functions are the same, as is shown in
great detail in \cit!{18}, 10.11--10.13 for the
analogous situation for smooth mappings.
\qed\enddemo

\proclaim{\nmb.{8.6}. Lemma (Curves in spaces of mappings)} Let
$M$ and $N$ be finite dimensional real analytic manifolds with
$M$ compact. \newline
1. A curve $c:\Bbb R\to C^\om(M,N)$ is real analytic if
and only if the associated mapping
$\hat c:\Bbb R\x M \to N$ is real analytic.

The curve $c:\Bbb R\to C^\om(M,N)$ is smooth if and only
if $\hat c:\Bbb R\x M\to N$ satisfies
the following condition:
\roster
\item"" For each $n$ there is an open neighborhood $U_n$ of
$\Bbb R\x M$ in $\Bbb R\x M_{\Bbb C}$ and a (unique)
$C^n$-extension $\tilde c: U_n\to N_{\Bbb C}$ such that 
$\tilde c(t,\quad)$ is holomorphic for all $t\in \Bbb R$.
\endroster
2. A curve $c:\Bbb R\to C^\infty(M,N)$ is smooth 
if and only if $\hat c:\Bbb R\x M\to N$
is smooth. 

The curve $c:\Bbb R\to C^\infty(M,N)$ is real analytic if and only
if $\hat c$ satisfies the following condition:
\roster
\item"" For each $n$ there is an open neighborhood $U_n$ of 
$\Bbb R\x M$ in $\Bbb C \x M$ and a (unique) $C^n$-extension $\tilde
c:U_n\to N_{\Bbb C}$ such that $\tilde c(\quad,x)$ is
holomorphic for all $x\in M$.
\endroster
\endproclaim

\demo{Proof} This follows from \nmb!{7.7} and the chart
structure on $C^\om(M,N)$. 
\qed\enddemo

\proclaim{\nmb.{8.7}. Corollary} Let $M$ and $N$ be real
analytic finite dimensional manifolds with $M$ compact. Let
$(U_\al,u_\al)$ be a real analytic atlas for $M$ and let $i:N\to
\Bbb R^n$ be a closed real analytic embedding into some $\Bbb R^n$.
Let $\Cal M$ be a possibly infinite dimensional real analytic manifold.

Then $f: \Cal M\to C^\om(M,N)$ is real analytic or smooth if and only
if 
$$C^\om(u_\al\i,i)\o f: \Cal M\to C^\om(u_\al(U_\al),\Bbb R^n)$$ 
is real analytic or smooth, respectively.

Furthermore $f: \Cal M\to C^\infty(M,N)$ is real analytic or smooth if 
and only if
$$C^\infty(u_\al\i,i)\o f: \Cal M\to C^\infty(u_\al(U_\al),\Bbb R^n)$$ 
is real analytic or smooth, respectively.
\endproclaim

\demo{Proof} By \nmb!{8.1} we may assume that $\Cal M =\Bbb R$.
Then we can use lemma \nmb!{8.6} on all appearing function spaces.
\qed\enddemo

\proclaim{\nmb.{8.8}. Theorem (Exponential law)} Let $\Cal M$ be
a (possibly infinite 
dimensional) real analytic manifold, and let $M$ and $N$ be
finite dimensional real analytic manifolds where $M$ is compact.

Then real analytic mappings $f:\Cal M\to C^\om(M,N)$ and real
analytic mappings $\widehat f:\Cal M\x M\to N$ correspond to
each other bijectively.
\endproclaim
\demo{Proof} Clearly we may assume that $\Cal M$ is a
$c^\infty$-open subset in a convenient vector space. Lemma
\nmb!{8.7} then reduces the assertion to cartesian closedness,
which holds by \nmb!{5.12}.
\qed\enddemo

\proclaim{\nmb.{8.9}. Corollary} If $M$ is compact and $M$, $N$ are
finite dimensional real analytic manifolds, then the evaluation
mapping $\ev: C^\om(M,N)\x M\to N$ is real analytic.

If $P$ is another compact real analytic manifold, then the
composition mapping
$\operatorname{comp}:C^\om(M,N)\x C^\om(P,M)\to C^\om(P,N)$ 
is real analytic.

In particular $f_*:C^\om(M,N)\to C^\om(M,N')$ and 
$g^*:C^\om(M,N)\to C^\om(P,N)$ are real analytic for 
$f\in C^\om(N,N')$ and $g\in C^\om(P,M)$.
\endproclaim
\demo{Proof} The mapping $\ev\spcheck = Id_{C^\om(M,N)}$ is real analytic,
so $\ev$ is it by \nmb!{8.8}. 
The mapping $\operatorname{comp}^\wedge = \ev\o(Id_{C^\om(M,N)}\x\ev):
C^\om(M,N)\x C^\om(P,M)\x P\to C^\om(M,N)\x M\to N$ is real analytic,
so also $\operatorname{comp}$.
\qed\enddemo

\proclaim{\nmb.{8.10}. Lemma} Let $M_i$ and $N_i$ are
finite dimensional real analytic manifolds with $M_i$ compact.
Then for $f\in C^\infty(N_1,N_2)$ the push forward
$f_*:C^\infty(M,N_1)\to C^\infty(M,N_2)$ is real analytic if and
only if $f$ is real analytic. 
For $f\in C^\infty(M_2,M_1)$ the pullback
$f^*:C^\infty(M_1,N)\to C^\infty(M_2,N)$ is, however, always
real analytic. 
\endproclaim 

\demo{Proof} If $f$ is real analytic and if $g\in C^\om(M,N_1)$,
then the mapping $u_{f\o g}\o
f_*\o u_g\i$  is a push forward by a real
analytic mapping, which is real analytic by \nmb!{7.8}. 

Obviously the canonical maps
$\operatorname{const}:N_1\to C^\infty(M,N_1)$ and
$\ev_x:C^\infty(M,N_2)\to N_2$ are real analytic. If $f_*$ is
real analytic, also $f = \ev_x\o f_*\o \operatorname{const}$ is it.

For the second statement choose real analytic atlases
$(U^i_\al,u^i_\al)$ of $M_i$ 
such that $f(U^2_\al)\subseteq U^1_\al$ and a closed real
analytic embedding $j:N\to \Bbb R^n$. Then the diagram
$$\define\toparrow{@>\pretend{f^*}
    \haswidth{(u^2_\al\o f\o (u^1_\al)\i)^*}>>}
\CD
C^\infty(M_1,N)  \toparrow    C^\infty(M_2,N)\\
@V{C^\infty((u^1_\al)\i,j)}VV    @VV{C^\infty((u^2_\al)\i,j)}V\\
C^\infty(u^1_\al(U^1_\al),\Bbb R^n)   
    @>>(u^2_\al\o f\o (u^1_\al)\i)^*>  C^\infty(u^2_\al(U^2_\al),\Bbb R^n)
\endCD$$
commutes, the bottom arrow is a continuous and linear mapping, so
it is real analytic. Thus by \nmb!{8.7} the mapping $f_*$ is
real analytic.
\qed\enddemo

\proclaim{\nmb.{8.11}. Theorem (Real analytic diffeomorphism
group)} For a compact real analytic manifold 
$M$ the group $\operatorname{Diff}^\om(M)$ of all real analytic
diffeomorphisms of $M$ is an open submani\-fold of $C^\om(M,M)$,
composition and inversion are real analytic.
\endproclaim
\demo{Proof} $\operatorname{Diff}^\om(M)$ is open in
$C^\om(M,M)$ in the compact open topology, thus also in the
finer manifold topology. The composition is real analytic by
\nmb!{8.9}, so it remains to show that the inversion 
$inv$ is real analytic. 

Let $c:\Bbb R\to
\operatorname{Diff}^\om(M)$ be a $C^\om$-curve. Then the
associated mapping $\hat c:\Bbb R\x M\to M$ is $C^\om$ by \nmb!{8.8} and 
$(inv\o c)^\wedge$ is the solution of the implicit equation 
$\hat c(t,(inv\o c)^\wedge (t,x)) = x$ and therefore real
analytic by the finite dimensional implicit function theorem.
Hence $inv\o c:\Bbb R\to \operatorname{Diff}^\om(M)$ is real analytic by
\nmb!{8.8} again.

Let $c:\Bbb R\to
\operatorname{Diff}^\om(M)$ be a $C^\infty$-curve. Then by
lemma \nmb!{8.6} the
associated mapping $\hat c:\Bbb R\x M\to M$ has a unique
extension to a $C^n$-mapping $\Bbb R\x M_{\Bbb C}\supseteq J\x
W\to M_{\Bbb C}$ which is holomorphic in the second variables
(has $\Bbb C$-linear derivatives),
for each $n\geq1$. The same assertion holds for the curve $inv\o
c$ by the finite dimensional implicit function theorem for $C^n$-mappings.
\qed\enddemo

\proclaim{\nmb.{8.12}. Theorem (Lie algebra of the diffeomorphism
group)} For a compact real analytic manifold 
$M$ the Lie algebra of the real analytic infinite dimensional
Lie group $\operatorname{Diff}^\om(M)$ is the convenient vector
space $C^\om(TM)$ of all real analytic vector fields on $M$,
equipped with the negative of the usual Lie bracket. The
exponential mapping $\operatorname{Exp}:C^\om(TM)\to
\operatorname{Diff}^\om(M)$ is the flow mapping to time 1, and
it is real analytic.
\endproclaim
\demo{Proof} The tangent space at $Id_M$ of
$\operatorname{Diff}^\om(M)$ is the space $C^\om(TM)$ of real
analytic vector fields on $M$, by \nmb!{8.5}. The one parameter
subgroup of a tangent vector is the flow $t\mapsto Fl^X_t$ of the
corresponding vector field $X\in C^\om(TM)$, so
$\operatorname{Exp}(X)=Fl^X_1$ which exists since $M$ is compact.

In order to show that $\operatorname{Exp}:C^\om(TM)\to
\operatorname{Diff}^\om(M)\subseteq C^\om(M,M)$ is real analytic,
by the exponential law \nmb!{8.8} it suffices to show that the
associated mapping 
$$\operatorname{Exp}^\wedge = Fl_1:C^\om(TM)\x M\to M$$
is real analytic. This follows from the finite
dimensional theory of ordinary real analytic and smooth
differential equations.  

For $X\in C^\om(TM)$ let $L_X$ denote the left invariant vector
field on $\operatorname{Diff}^\om(M)$. Its flow is given by 
$Fl^{L_X}_t(f)=f\o\operatorname{Exp}(tX)$.
The usual proof of differential geometry shows that 
$[L_X,L_Y]=\tfrac d{dt}|_0 (Fl^{L_X}_t)^*L_Y$, thus for $e=Id_M$
we have
$$\allowdisplaybreaks\align 
[L_X,L_Y](e)&= (\tfrac d{dt}|_0 (Fl^{L_X}_t)^*L_Y)(e)\\
&= \tfrac d{dt}|_0(T(Fl^{L_X}_{-t})\o L_Y\o Fl^{L_X}_t)(e)\\
&= \tfrac d{dt}|_0T(Fl^{L_X}_{-t})(L_Y(e\o Fl^X_t)\\
&= \tfrac d{dt}|_0T(Fl^{L_X}_{-t})(T(Fl^X_t)\o Y)\\
&= \tfrac d{dt}|_0(T(Fl^X_t)\o Y\o Fl^X_{-t})\\
&= \tfrac d{dt}|_0(Fl^X_{-t})^*Y = -[X,Y].\quad\qed
\endalign$$
\enddemo

\proclaim{\nmb.{8.13}. Example} The exponential map
$\operatorname{Exp}:C^\om(TS^1)\to \operatorname{Diff}^\om(S^1)$
is neither locally injective nor surjective on any neighborhood
of the identity.   
\endproclaim
\demo{Proof} The proof of \cit!{24}, 3.3.1
for the group of smooth diffeomorphisms of $S^1$ can be adapted
to the real analytic case:
$$\ph_n(\th) = \th + \tfrac{2\pi}n + \tfrac1{2^n}\sin n\th$$
is Mackey convergent (in $U_{Id}$) to $Id_{S^1}$ in
$\operatorname{Diff}^\om(S^1)$ and is not in the image of the
exponential mapping. 
\qed\enddemo  

\subheading{\nmb.{8.14}. Remarks}
For a real analytic manifold $M$ the group
$\operatorname{Diff}(M)$ of all smooth diffeomorphisms of $M$ is
a real analytic open submanifold of $C^\infty(M,M)$ and is a
smooth Lie group by \cit!{18}, 11.11. The composition mapping is
not real analytic by \nmb!{8.10}. Moreover it does not
carry any real analytic Lie group structure by
\cit!{22}, 9.2, and it has no complexification in
general, see \cit!{24}, 3.3.  The mapping
$$Ad\o\operatorname{Exp}:C^\infty(TM)\to \operatorname{Diff}(M)\to
L(C^\infty(TM),C^\infty(TM))$$
is not real analytic, see
\cit!{19}, 4.11.

For $x\in M$ the mapping 
$ev_x\o \operatorname{Exp}:C^\infty(TM)\to
\operatorname{Diff}(M)\to M$ is not real analytic, since 
$(ev_x\o\operatorname{Exp})(tX) = Fl^X_t(x)$ which is not real
analytic in $t$ for general smooth $X$.

The exponential mapping $\operatorname{Exp}:C^\infty(TM)\to
\operatorname{Diff}(M)$ is in a very strong
sense not surjective: In \cit!{6} it is shown, that
$\operatorname{Diff}(M)$ contains an arcwise connected free
subgroup on $2^{\aleph_0}$ generators which meets the image of
$\operatorname{Exp}$ only at the identity. 

The real analytic Lie group $\operatorname{Diff}^\om(M)$ is \idx{\it
regular} in the sense of \cit!{22}, 7.6, where
the original concept of \cit!{23} is weakened. 
This condition means
that the mapping associating the evolution operator to each time
dependent vector field on $M$ is smooth. It is even real
analytic, compare the proof of theorem \nmb!{8.12}.

\proclaim{\nmb.{8.15}. Theorem (Principal bundle of embeddings)}
Let $M$ and $N$ be real 
analytic manifolds with $M$ compact. Then the set $Emb^\om(M,N)$
of all real analytic embeddings $M\to N$ is an open submanifold
of $C^\om(M,N)$. It is the total space of a real analytic
principal fiber bundle with structure group
$\operatorname{Diff}^\om(M)$, whose real analytic base manifold
is the space of all submanifolds of $N$ of type $M$.
\endproclaim
\demo{Proof} The proof given in \cit!{18}, section 13
or \cit!{5}, 4.7.8
is valid with the obvious changes. One starts with a real
analytic Riemannian metric and uses its exponential mapping.
The space of embeddings is open, since embeddings are open in
$C^\infty(M,N)$, which induces a coarser topology.
\qed\enddemo

\proclaim{\nmb.{8.16}. Theorem (Classifying space for
$\operatorname{Diff}^\om(M)$)}  Let $M$ be a compact real
analytic manifold. Then the space $Emb^\om(M,\ell^2)$ of
real analytic embeddings of $M$ into the
Hilbert space $\ell^2$ is the total space of a real analytic principal fibre
bundle with structure group $\operatorname{Diff}^\om(M)$ and
real analytic base manifold $B(M,\ell^2)$, which is a classifying
space for the Lie group $\operatorname{Diff}^\om(M)$.
\endproclaim

\demo{Proof} The construction in \nmb!{8.15}
carries over to the Hilbert space $N=\ell^2$ with the appropriate
changes to obtain a real analytic principal fibre bundle. Its
total space is continuously contractible and so the bundle is
classifying, see the argument in \cit!{20}, section 6.
\qed\enddemo


\Refs
\widestnumber\key{23}
\ref
\key\cit0{1}
\by Abraham, R.  \book Lectures of Smale on differential
topology \bookinfo Lecture Notes \publ Columbia University
\publaddr New York \yr 1962 \endref

\ref 
\key\cit0{2}
\by Bochnak, J.\& Siciak, J. 
\paper Analytic functions in
topological vector spaces \jour Studia Math. \vol 39 \yr 1971
\pages 77--112\endref

\ref 
\key\cit0{3}
\by Boman, J. 
\paper Differentiability of a function and
of its compositions with functions of one variable \jour Math.
Scand. \vol 20 \yr 1967 \pages 249--268 \endref

\ref 
\key\cit0{4}
\by Floret, K. 
\paper Lokalkonvexe Sequenzen mit
kompakten Abbildungen \jour J. Reine Ang. Math. \vol 247 \yr
1971 \pages 155--195 \endref

\ref 
\key\cit0{5} 
\by Fr\"olicher, A. \& Kriegl, A. 
\book Linear spaces and differentiation theory 
\bookinfo Pure and Applied Mathematics 
\publ J. Wiley \publaddr Chichester \yr 1988 \endref

\ref
\key\cit0{6}
\by Grabowski, J. \paper Free subgroups of
diffeomorphism groups  \jour Fundamenta Math. \vol 131 
\pages 103--121 \yr 1988 \endref

\ref 
\key\cit0{7}
\by Grauert, H. 
\paper On Levi's problem and the embedding
of real analytic manifolds \jour Annals of Math. \vol 68 \pages
460--472 \yr 1958 \endref

\ref 
\key\cit0{8}
\by Grothendieck, A. 
\paper Sur certains espaces de fonctions holomorphes. I 
\jour J. Reine Ang. Math. \vol 192 
\yr 1953 \pages 35--64 \endref

\ref 
\key\cit0{9}
\by H\"ormander, L. \book The Analysis of linear partial
differential operators I \bookinfo Grund\-lehren 256 \publ
Springer-Verlag \publaddr Berlin \yr 1983 \endref

\ref 
\key\cit0{10} 
\by Jarchow, H. \book Locally convex spaces
\publ Teubner \publaddr Stuttgart \yr 1981 \endref

\ref 
\key\cit0{11}
\by Kashiwara, M., Kawai, T. \& Kimura, T.
\book Foundations of algebraic analysis \publ Princeton Univ.
Press \publaddr Princeton \yr 1986 \endref

\ref 
\key\cit0{12}
\by K\"othe, G. \paper Dualit\"at in der
Funktionentheorie  \jour J. Reine Ang. Math. \vol 191 \yr 1953
\pages 30--49 \endref

\ref 
\key\cit0{13}
\by Kriegl, A. \paper Die richtigen R\"aume f\"ur Analysis
im Unendlich - Dimensionalen \jour Monatshefte Math. \vol 94 \yr
1982 \pages 109--124 \endref

\ref 
\key\cit0{14}
\by Kriegl, A. \paper Eine kartesisch abgeschlossene
Kategorie glatter Abbildungen zwischen beliebigen lokalkonvexen
Vektorr\"aumen \jour Monatshefte f\"ur Math. \vol 95 \yr 1983 \pages
287--309 \endref

\ref
\key\cit0{15} 
\by Kriegl, A. \& Nel, L. D. \paper A convenient
setting for holomorphy \jour Cahiers Top. G\'eo. Diff. \vol 26
\yr 1985 \pages 273--309 \endref

\ref 
\key\cit0{16}
\by Leslie, J. \paper On the group of real analytic
diffeomorphisms of a compact real analytic manifold \jour
Transact. Amer. Math. Soc. \vol 274 (2) \yr 1982 \pages 651--669 \endref

\ref
\key\cit0{17}
\by Mazur, S. \& Orlicz, W. \paper Grundlegende Eigenschaften
der polynomischen Operationen \jour Studia Math. \vol 5 \yr 1935
\pages 50--68 \endref

\ref
\key\cit0{18}
\by Michor, P. W. \book Manifolds of
differentiable mappings \publ Shiva \yr 1980 \publaddr Orpington \endref

\ref 
\key\cit0{19}
\by Michor, P. W. \paper Manifolds of smooth
mappings IV: Theorem of De~Rham \jour Cahiers Top. Geo. Diff. 
\vol 24 \yr 1983 \pages 57--86 \endref

\ref 
\key\cit0{20}
\by Michor, P. W. \paper Gauge theory for the
diffeomorphism group \inbook Proc. Conf.
Differential Geometric Methods in Theoretical Physics, 
Como 1987, K. Bleuler and M. Werner (eds)\publ Kluwer
\publaddr Dordrecht  \yr 1988 \pages 345--371 \endref

\ref
\key\cit0{21}
\by Michor, P. W. \paper The moment mapping for unitary
representations \paperinfo Preprint 1989 \endref

\ref
\key\cit0{22}
\by Milnor, J. \paper Remarks on infinite dimensional Lie groups
\inbook Relativity, Groups, and Topology II, Les Houches, 1983,
B.S.~DeWitt, R.~Stora, Eds. \publ Elsevier \yr 1984 \publaddr
Amsterdam \endref 

\ref 
\key\cit0{23}
\by Omori, H., Maeda, Y., Yoshioka, A. \& Kobayashi, O.
\paper On regular Fr\'echet Lie groups IV, Definition and fundamental theorems
\jour Tokyo J. Math.
\vol 6
\pages 39--64
\yr 1988
\endref

\ref 
\key\cit0{24}
\by Pressley, A. \& Segal, G. 
\book Loop groups 
\bookinfo Oxford Mathematical Monographs \publ Oxford University
Press \yr 1986 \endref

\ref 
\key\cit0{25}
\by Sebasti\~ao e Silva, J. \paper As fun\c c\~oes
anal\'iticas e a an\'alise funcional \jour Port. Math. \vol 9
\yr 1950 \pages 1--130 \endref

\ref 
\key\cit0{26}
\by Sebasti\~ao e Silva, J. \paper Sobre a topologia
dos espa\c cos funcionais anal\'iticos \jour Rev. Fac. Cienc.
Lisboa \vol 2, ser. A1 \yr 1950 \pages 23--102 \endref

\ref 
\key\cit0{27}
\by Sebasti\~ao e Silva, J. 
\paper Su certe classi di spazi localmente convessi importanti per le 
applicationi 
\jour Rendiconti di Matematica, Roma, Ser. V
\vol 14
\pages 388--410
\yr 1955
\endref

\ref 
\key\cit0{28}
\by Semmes, S. \paper Nonlinear Fourier analysis \jour
Bull.(New Series) Amer. Math. Soc. \vol 20 \pages 1--18 \yr 1989
\endref

\ref 
\key\cit0{29}
\by Toeplitz, O. \paper Die linearen vollkommenen R\"aume
der Funktionentheorie \jour Comm. Math. Helv. \vol 23 \yr 1949 \pages
222--242 \endref

\ref 
\key\cit0{30}
\by Treves, F. \book Topological vector spaces,
distributions, and kernels \publ Academic Press \publaddr New
York \yr 1967 \endref

\ref 
\key\cit0{31}
\by Van Hove, L. \paper Topologie des espaces
fontionnels analytiques et de groupes infinis de transformations
\jour Bull. Classe des Sciences, Acad. Royale Belgique \vol 38
\yr 1952 \pages 333--351 \endref

\ref 
\key\cit0{32}
\by Whitney, H. \& Bruhat, F. \paper Quelques
propri\'et\'es fondamentales des ensembles analytiques-r\'eels
\jour Comm. Math. Helv. \vol 33  \yr 1959 \pages 132--160 \endref

\endRefs
\enddocument